\documentclass[12pt,reqno]{amsart}
\usepackage{amsmath}
\usepackage{amssymb}
\usepackage{amsthm}
\usepackage{mathrsfs}
\usepackage{latexsym}
\usepackage{amssymb}

\newtheorem{theorem}{Theorem}[section]
\newtheorem{lemma}{Lemma}[section]
\newtheorem{corollary}{Corollary}[section]
\newtheorem{remark}{Remark}[section]
\newtheorem{definition}{Definition}[section]
\newtheorem{proposition}{Proposition}[section]
\newtheorem{example}{Example}[section]
\newtheorem{assumption}{Assumption}[section]
\numberwithin{equation}{section}
\newcommand{\bth}{\begin{theorem}}
	\newcommand{\ethe}{\end{theorem}}

\newcommand{\bre}{\begin{remark}}
	\newcommand{\ere}{\end{remark}}

\newcommand{\ble}{\begin{lemma}}
	\newcommand{\ele}{\end{lemma}}
\newcommand{\bde}{\begin{definition}}
	\newcommand{\ede}{\end{definition}}

\newcommand{\bco}{\begin{corollary}}
	\newcommand{\eco}{\end{corollary}}

\newcommand{\bpr}{\begin{proposition}}
	\newcommand{\epr}{\end{proposition}}

\newcommand{\bexer}{\begin{exercise}}
	\newcommand{\eexer}{\end{exercise}}
\newcommand{\breh}{\begin{hint}}
	\newcommand{\ereh}{\end{hint}}

\newcommand{\halmos}{\hfill \qed}

\newcommand{\bexam}{\begin{example}}
	\newcommand{\eexam}{\end{example}}

\newcommand{\pr} {{\bf Proof.}}

\newcommand{\bfi}{\begin{fig}}
	\newcommand{\efi}{\end{fig}}

\newcommand{\beao}{\begin{eqnarray*}}
	\newcommand{\eeao}{\end{eqnarray*}\noindent}

\newcommand{\beam}{\begin{eqnarray}}
	\newcommand{\eeam}{\end{eqnarray}\noindent}
\newcommand{\E}{\mathbf{E}}
\newcommand{\PP}{\mathbf{P}}

\newcommand{\xto}{x\to\infty}

\newcommand{\bF}{\overline{F}}

\newcommand{\bG}{\overline{G}}
\newcommand{\bV}{\overline{V}}

\newcommand{\bbr}{{\mathbb R}}

\newcommand{\bbb}{{\mathbb B}}

\newcommand{\bbn}{{\mathbb N}}

\newcommand{\vep}{\varepsilon}

\usepackage{soul}
\setstcolor{red}

\usepackage{color}

\begin{document}
	\title[The Bivariate Regular Variation weighted sums]{The Bivariate Regular Variation of randomly weighted sums revisited in the presence of interdependence}

\author[ D.G. Konstantinides, C. D. Passalidis ]{ Dimitrios G. Konstantinides, Charalampos  D. Passalidis} 

\address{Dept. of Statistics and Actuarial-Financial Mathematics,
University of the Aegean,
Karlovassi, GR-83 200 Samos, Greece}
\email{konstant@aegean.gr,\;sasd24009@sas.aegean.gr.}
	\date{{\small \today}}
	
\begin{abstract} 
We study the joint distribution of two randomly weighted sums. Inspired by the practical applications, we assume that the main random variables follow the non-standard bivariate regular variation, symbolically $BRV$, to put emphasis to the value of inhomogeneity of the risk distribution tails, while the random weights are weakly dependent with main random variables. Under some moment conditions on the random weights we show that the randomly weighted sums have $BRV$ distribution with an analytic relation for the Radon measure, that captures the interdependence between the random weights and the main random variables. Under some stronger moments conditions, our result is extended uniformly, with respect to summands, covering also the case of infinite randomly weighted sums. In order to keep weak dependence structure among the random weights and the main random variables, we require the random weights to be independent each other, something that does not happen in models with insurance and financial risks. Up to recent years, such kind of approximations, even in one-dimensional case, had mostly theoretical interest, since underline the presence of (multivariate linear) single big jump principle. However, here we provide an application of the main results on ruin probability in a new flexible credit risk model. In our model, although we restrict ourselves to standard $BRV$, the obliged do not enter - quit necessarily simultaneously to the system, while the breach probability is not necessarily independent of the the amount of breach for the obliged. Finally, in the non-standard $BRV$, with asymptotically dependent risks, we provide an application of the main results, to find the asymptotic behavior of a risk measure, which is called joint expected shortfall, that plays crucial role to the measure of the contagion of extreme risks. 
\end{abstract}
	
	\maketitle
	\vspace{3mm}
	

\section{Introduction} \label{sec.KLP.1}

	Let $\{(X_i,\,Y_i)\,,\;i \in \bbn\}$ be a sequence of independent, identically distributed (i.i.d.) random vectors, which are coppies of the general random vector $(X,\,Y)$, with distributions $(F,\,G)$.  We assume that the $(X,\,Y)$ represents a pair of nonnegative random variables with tails $\bF=1-F$, $\bG=1-G$ such that the inequalities $\bF(x)>0$, $\bG(x)>0$ are true for any $x \in \bbr:=(-\infty,\,\infty)$. Let us also consider
	the sequences $\{\Theta_i, \, \Delta_i,\; i \in \bbn\}$ of arbitrarily dependent, nonnegative and non-degenerate at zero random variables. We assume that the vector $(X,\,Y)$ follows a bivariate regular varying ($BRV$) distribution, see in Subsection $2.1$ below, and we are interested in the behavior of the pair $(S_n^{\Theta},\,T_n^{\Delta})$, in the sense of $BRV$ structure, where 
	\beam \label{eq.AKP.1.1}
	S_{n}^\Theta:=\sum_{i=1}^{n}\Theta_i X_i\,,\qquad T_{n}^\Delta:=\sum_{i=1}^{n}\Delta_i Y_i\,,
	\eeam
	for $n \in \bbn$ (as also for $n=\infty$), under some dependence conditions on $(X_i,\,\Theta_i)$ and on $(Y_i,\,\Delta_i)$.
	
	Although the asymptotic behavior of the tail of a randomly weighted sum $S_{n}^\Theta$ was studied extensively, see for example \cite{gao:wang:2010},\cite{yang:leipus:siaulys:2012},\cite{li:2013},\cite{yang:wang:leipus:siaulys:2013},\cite{cheng:2014}, among others, only few papers spend attention on the two dimensional case. Even more in these few publications the independence assumption among the main variables $\{(X_i,\,Y_i)\,,\;i \in \bbn\}$  and among the discount factors (or, random weights) $\{\Theta_i,\,\Delta_i,\; i \in \bbn\}$ makes the condition of practical usage to be pretty narrow as in case of independence between insurance and financial risks, see for example \cite{li:2018}, \cite{chen:yang:2019}, \cite{shen:du:2023}, \cite{yang:chen:yuen:2024}, \cite{chen:cheng:zheng:2025}, \cite{konstantinides:passalidis:2024c}.   
	
To depict the full picture in the current literature in (as one-dimensional, as also multidimensional) randomly weighted sums, there exist two main directions.
\begin{enumerate}
\item[(i)]
The $\{(X_i,\,Y_i),\,i \in \bbn\}$ have some dependence structure, as also the $\{(\Theta_i,\,\Delta_i),\,i \in \bbn\}$ (usually arbitrarily dependent), but the sequence $\{(X_i,\,Y_i),\,i \in \bbn\}$ is independent of  $\{(\Theta_i,\,\Delta_i),\,i \in \bbn\}$, see for example \cite{li:2018}, \cite{chen:yang:2019}, \cite{konstantinides:passalidis:2024c}. 
\item[(ii)]
The  $\{(X_i,\,Y_i,\,\Theta_i,\,\Delta_i),\; i \in \bbn\}$ have some dependence structure, but this sequence has independent terms.
\end{enumerate}

With respect to direction (ii), we know only a few contributions on one-dimensional set up, see for example \cite{yang:leipus:siaulys:2012}, \cite{yang:wang:leipus:siaulys:2013}. Additionally, these contributions possess mostly theoretical motivation, but the independence of random weights $\{(\Theta_i,\,\Delta_i),\,i \in \bbn\}$ makes then practically problematic in applications of usual risk models (discrete-time risk models with insurance and financial risks). As was mentioned in \cite[Sec. 5]{chen:cheng:zheng:2025}, the direction (ii) in two-dimensional set up is very challenging.

In this paper we focus in direction (ii), in fact complementing direction (i) for $BRV$, that was given in \cite{chen:yang:2019}. Furthermore, we provide an application in ruin probability of an insurer, who faces credit risk. Such an application gives the essence of the contribution of direction (ii) in actuarial science.
  
We feel as necessary to take into account the essential dependence between the primary random variables and random weights from direction (ii) (as we shall show in Sec. $6$), from where we get motivated to examine the closure properties of $BRV$ for the pair $(S_n^{\Theta},\,T_n^{\Delta})$, with $n \in \bbn \cup \{\infty\}$, under a weak dependence structure among the random variables $X_i$, $\Theta_i$, $Y_i$, $\Delta_i$,  that contains independence as a special case, and also allows the vectors $(X_i,Y_i)$ and $(\Theta_{i}\,\, \Delta_{i})$ to have arbitrarily dependent components. We should recall that the simultaneous dependence between the main random variables $(X,\,Y)$, but also the dependence of them with the random weights $(\Theta,\, \Delta)$ is called interdependence, as it appears in \cite{konstantinides:leipus:passalidis:siaulys:2024}, where the one-dimensional asymptotic behavior of $S_n^{\Theta}$ was studied, see also \cite{chen:cheng:2024} for similar studies on interdependence.  
	
	Nevertheless,  we meet very often the interdependence effect in multivariate risk models and in risk measures with dependent portfolios, see for example \cite{fougeres:mercadier:2012}, \cite{cheng:konstantinides:wang:2024}, \cite{konstantinides:passalidis:2024a}. However, in all these papers was considered only the standard case of multivariate regular variation.
	
	As main domain of application of the closure properties about the pair $(S_n^{\Theta},\,T_n^{\Delta})$ with condition of $BRV$ and under some interdependence structures can be the risk theory and risk management and our applications of main results, also focus on these two topics. For some results that focus on the closure properties for heavy tailed distributions we refer the reader to \cite{leipus:siaulys:konstantinides:2023} for the one-dimensional case and \cite{fougeres:mercadier:2012}, \cite{samorodnitsky:sun:2016}, \cite{das:fashenhartmann:2023}, \cite{konstantinides:passalidis:2024b}, \cite{konstantinides:passalidis:2024h} for the multivariate ones.
	
	The rest of this paper is organized as follows. In Section 2, we provide the necessary preliminary results either for the distribution classes or for the dependence structures. In Section 3, we present the first main result, as an  extension of Breiman's theorem in $BRV$ set up under interdependence. In Section 4, we establish the closure property of distribution class $BRV$ with respect to pair $(S_n^{\Theta},\,T_n^{\Delta})\,,\;n \in \bbn$. In Section $5$, we present the closure properties of $BRV$ with respect to $(S_n,\,T_n)$, for $n \in \bbn \cup \{\infty\}$, uniformly for all $n \in \bbn$, assuming stronger moment conditions, than that from Section $4$. In Section $6$ are presented two applications of our main results. The first one relates to the Risk Theory. More precisely, restricting ourselves in standard $BRV$ case, was established the asymptotic behavior of three types of ruin probability, in a credit risk model, under the presence of interdependence. In this model, the potential obligators do not enter or quit necessarily at the same moments from the system, while the probability of breach and the amount of breach are weakly dependent. According to the best of our knowledge, such an approach of credit risk is new in the literature.  The second application relates to risk management. Under the $BRV$ structure (either standard or non-standard), when the products of these risks satisfy asymptotic dependence conditions, we provide the asymptotic behavior of a risk measure, called joint expected shortfall.

\section{Preliminaries} \label{sec.KLP.2}

All the limit relations hold as $\xto$, except stated otherwise. For any two positive functions $f,\,g$, the asymptotic relation 
\beao
f(x)=o[g(x)]\,,
\eeao 
means that we have
\beao
\lim \dfrac{f(x)}{g(x)} = 0\,.
\eeao
The asymptotic relation  $f(x)=O[g(x)]$, means that we have
\beao
\limsup \dfrac{f(x)}{g(x)} < \infty\,,
\eeao
while the asymptotic relation  $f(x)\asymp g(x)$, means that both $f(x)=O[g(x)]$, and $g(x)=O[f(x)]$ hold. For any event $E$, we denote by ${\bf 1}_{E}$ the indicator function of the event $E$. For a set $E$,  we denote by $E^c$ its complement and by $\overline{E}$ its closed hull. For any $x,\,y$ real numbers, we denote $x\vee y :=\max\{x,\,y\}$ and $x \wedge y := \min\{x,\,y\}$, $x^+:=x \vee 0$. With ${\bf 0}$ we denote the zero vector. For two bivariate vectors ${\bf a}=(a_1,\,a_2)^{\top}$, ${\bf b}=(b_1,\,b_2)^{\top}$ we denote by ${\bf a} \leq {\bf b}$ (or, equivalently, ${\bf b} \geq {\bf a}$), having in mind component-wise comparison $a_i \leq b_i$, for $i=1,\,2$.
	
The random variables follow distributions with support on $\bbr_+:=[0,\,\infty)$, unless otherwise stated. The support of a distribution $V$ is denoted by $S_V$ and the random variable $Z$ has distribution $V$.

\subsection{Bivariate regular variation} \label{sec.KP.2.1}

For a distribution $V$, such that $\bV(x)>0$ for any $x \in \bbr$, we say that it is regularly varying, with some finite index $\alpha>0$ if it holds
\beam \label{eq.AKP.2.2}
\lim \dfrac{\bV(t\,x)}{\bV(x)}=t^{-\alpha}\,,
\eeam
for any $t>0$, and we write symbolically $V \in \mathcal{R}_{-\alpha}$. The class of regular variation plays crucial role in the frame of heavy tailed distribution, with several application in actuarial science, finance, credit risk and risk management, see for example \cite{konstantinides:mikosch:2005}, \cite{zhu:li:2012}, \cite{li:2022a}, \cite{li:2023}, \cite{liu:yi:2025} among others. For further treatments on regular variation, see \cite{bingham:goldie:teugels:1987} and \cite{resnick:1987}.
	
An equivalent definition of the regular variation, as it appears in \eqref{eq.AKP.2.2} is the following: There exists finite, non-degenerate to zero, Radon measure $\mu$, such that  the convergence
\beam \label{eq.AKP.2.3}
\lim_{\xto} x\,\PP\left[\dfrac Z{U_V(x)} \in \bbb \right]=\mu(\bbb)\,,
\eeam
holds for any Borel set $\bbb \subset \overline{\bbr}_+=[0,\,\infty]$, with ${\bf 0} \notin \overline{\bbb}$, which is $\mu$-continuous, where
\beao
U_V(x):=\left(\dfrac 1{\bV} \right)^{\leftarrow}(x)=\inf \left\{y \in \bbr \,:\,\dfrac 1{\bV(y)} \geq x \right\}\,,
\eeao
for any $x>0$. We call $U_V(x)$ normalization function. The measure $\mu$ possesses the property of homogeneity, namely for any $\lambda>0$ and any Borel set  $\bbb \subset \overline{\bbr}_+$, with ${\bf 0} \notin \overline{\bbb}$, it holds
\beao
\mu\left(\lambda^{1/\alpha} \bbb \right)=\dfrac 1{\lambda} \mu(\bbb)\,.
\eeao
	
The bivariate regular variation $BRV$, represents a special case of multivariate regular variation $MRV$, which introduced in \cite{dehaan:resnick:1981}. We say that a nonnegative random vector $(X,\,Y)$ follows a $BRV$ structure, if there exists a finite, nondegenerate to zero, Radon measure $\mu$, and two distributions $F$ and $G$, such that provides
\beam \label{eq.AKP.2.6}
\lim x\,\PP\left[\left(\dfrac X{U_F(x)},\;\dfrac Y{U_G(x)}  \right)\in \bbb \right]=\mu\left(\bbb \right)\,,
\eeam
for any Borel set  $\bbb \subset \overline{\bbr}_+^2:=[0,\,\infty]^2$, with ${\bf 0} \notin \overline{\bbb}$, which is $\mu$-continuous. The distributions $F$, $G$, appearing in the normalization functions, are such that $F \in \mathcal{R}_{-\alpha}$ and $G \in \mathcal{R}_{-\beta}$, for some finite indexes $\alpha, \,\beta >0$. For the distributions that satisfy \eqref{eq.AKP.2.6} we write symbolically $(X,\,Y) \in BRV_{-\alpha,\,-\beta}(\mu)$. We notice that if $U_V(x)=U_F(x)=U_G(x)$, that means
\beam \label{eq.AKP.2.7}
\lim  x\,\PP\left[\dfrac {(X,\,Y)}{U_V(x)}\in \bbb \right]=\mu\left(\bbb\right)\,,
\eeam
with $V \in \mathcal{R}_{-\alpha}$, then we obtain the standard $BRV$ structure, symbolically $(X,\,Y)\in BRV_{-\alpha}(\mu)$. Hence, the general $BRV$ structure in \eqref{eq.AKP.2.6} contains the standard $BRV$ structure, as special case. Here again we observe, that the measure $\mu$, possesses the property of homogeneity, in the sense that for any $\lambda>0$ and Borel set  $\bbb \subsetneq \overline{\bbr}_+^2$, with ${\bf 0} \notin \overline{\bbb}$, we have
\beao	
\mu (\bbb_{\lambda})=\dfrac 1{\lambda}\,\mu\left(\bbb \right)\,,	
\eeao
where $\bbb_{\lambda}=\{(\lambda^{1/\alpha}\,p,\;\lambda^{1/\beta}\,q)\;:\;(p,\,q)\in \bbb\}$, see also in \cite[Lem. 5.1]{tang:yang:2019}.
	
Several papers  with applications of $MRV$ on risk theory and risk management, adopt as condition the standard $MRV$, see for example \cite{asimit:furman:tang:vernic:2011}, \cite{li:2016}, \cite{li:2022b}, \cite{zhu:li:yang:xie:sun:2023}, \cite{chen:yang:zhang:2024}, \cite{chen:tong:yang:2023}, among others. However, there still exist a few papers with the condition of  general $MRV$, usually restricted in the bivariate case $BRV$, see \cite{chen:yang:2019}, \cite{tang:yang:2019}, \cite{tang:xun:zhou:2026}.
	
	The fact that the $BRV$ structure for the vector $(X,\,Y)$ permits arbitrary dependence between the pair $(X,\,Y)$, in combination with the possible variety in the tail indexes of the distributions of $X$ and $Y$, are the main reasons that $BRV$ attracted attention. Furthermore, we recognize these reasons, as motivation of the study of the closure properties in this distribution class, with hopeful impact in practical applications.

\subsection{The dependence concept} \label{subsec.KP.2.2}

Now, we give the main dependencies for further usage. Each pair $(X_i,\,Y_i)$ contains arbitrary dependence under the $BRV$ structure. Therefore, our interest is focused in the dependencies between $\Theta_i$ and $X_i$ as also between $\Delta_i$ and $Y_i$, in order to reach a generalization up to the interdependence.
	
The first dependence structure was introduced in \cite[Assum. 2]{asimit:jones:2008} and has be used in risk theory, to depict either time-dependent risk models, in continuous time risk models see \cite{asimit:badescu:2010},\cite{li:tang:wu:2010}, or to present dependence structure between the random weights and the main random sizes, in discrete time risk models, see for example \cite{yang:leipus:siaulys:2012},\cite{yang:wang:leipus:siaulys:2013}.
	
Let us render this dependence in the following assumption for any pair $(X_i,\,\Theta_i)$ and any pair $(Y_i,\,\Delta_i)$. This dependence is formulated for only one random pair, but we use later this kind of dependence on the sequences of such pairs. 
	
\begin{assumption} \label{ass.AKP.2.1}
\begin{enumerate}	
		
(a)
For any pair $(X_i,\,\Theta_i)$, with $i \in \bbn$ there exists some measurable function $h_{1,i}\;:\;[0,\,\infty) \to (0,\,\infty)$, such that
\beam \label{eq.AKP.2.9}
\PP[X_i>x\;\big|\;\Theta_i =\theta]\sim h_{1,i}(\theta)\,\PP[X_i>x]\,,
\eeam
uniformly for any $\theta \in S_{\Theta_i}$.
			
(b)
For any pair $(Y_i,\,\Delta_i)$, with $i \in \bbn$ there exists some measurable function $h_{2,i}\;:\;[0,\,\infty) \to (0,\,\infty)$, such that
\beam \label{eq.AKP.2.10}
\PP[Y_i>x\;\big|\;\Delta_i =\delta] \sim h_{2,i}(\delta)\,\PP[Y_i>x]\,,
\eeam
uniformly for any $\delta \in S_{\Delta_i}$.
\end{enumerate}
\end{assumption} 
	
Note that the uniformity property in \eqref{eq.AKP.2.9}  is understood as 
\beao
\lim \sup_{\theta \in S_{\Theta_i}} \left| \dfrac{\PP[X_i > x\;\big|\; \Theta_i=\theta]}{h_{1,i}(\theta)\,\PP[X_i>x]} -1 \right|=0\,,	
\eeao
and similarly in \eqref{eq.AKP.2.10}.
	
\bre \label{rem.AKP.2.1}
The notation $\Theta_i = \theta$ and $\Delta_i = \delta$, are technical, in the sense that the $\theta$ and $\delta$ are not necessarily points, but also can be intervals, in the sense that the left member in \eqref{eq.AKP.2.9} is defined by the relation
\beao
\lim_{\vep \downarrow 0} \PP[X_i > x\;\big|\;\Theta_i \in [\theta - \vep,\,\theta + \vep]]\,,
\eeao 
and similarly for relation \eqref{eq.AKP.2.10}. Furthermore, we can see that relations \eqref{eq.AKP.2.9} and \eqref{eq.AKP.2.10} imply as special case the independence between $(X_i,\,\Theta_i)$ and between $(Y_i,\,\Delta_i)$, respectively (in this case is implied $h_{1,i}\equiv h_{2,i} =1$, but the opposite does not gives necessarily the independence, see \cite[Prop. 2.6]{cui:wang:2025}). Additionally, by \cite[Prop. 2.4]{cui:wang:2025}, we get that the functions $h_{1,i}$ and  $h_{2,i}$ are bounded from above by some constant, namely there exist $k_{1,i}$, $k_{2,i}$, with $0<k_{1,i},\,k_{2,i} < \infty$, such that 
\beam \label{eq.AKP.2.11}
h_{1,i}(\theta_i)\leq k_{1,i}\,, \qquad h_{2,i}(\delta_i)\leq k_{2,i}\,,
\eeam
for any  $i \in \bbn$, with $\theta_i\in S_{\Theta_i}$ and $\delta_i\in S_{\Delta_i}$. Finally, it is easy to  see that
\beam \label{eq.AKP.2.12}
\E[h_{1,i}(\Theta_i)]=1=\E[h_{2,i}(\Delta_i)]\,,
\eeam
for any  $i \in \bbn$.
\ere
	
In the following statement, when say that holds Assumption \ref{ass.AKP.2.1} (a) for the pair  $(X,\,\Theta)$, it means that there exists a function $h_1\;:\;[0,\infty) \to (0,\,\infty)$ such that 
\beam \label{eq.AKP.2.1*}
\PP[X>x\;\big|\;\Theta =\theta]\sim h_{1}(\theta)\,\PP[X>x]\,,
\eeam
uniformly for any $\theta \in S_{\Theta}$.
	
In the following proposition, we find that an extension of the Breiman's theorem, under the dependence structure Assumption \ref{ass.AKP.2.1} (a), found in \cite[Th. 2.1 (i)]{cui:wang:2025}. For the original version of Breiman's theorem, when the $X,\,\Theta$ are independent, see \cite{breiman:1965},\cite{cline:samorodnitsky:1994} and in generalization with dependencies, see \cite{yang:wang:2013},\cite{cui:wang:2025}.
	
\bpr \label{pr.AKP.2.1}
Let $X$ be a real-valued random variable with distribution $F \in \mathcal{R}_{-\alpha}$, with $\alpha>0$ and $\Theta$ be nonnegative, nondegenerated to zero random variable. If the pair $(X,\,\Theta)$ satisfies Assumption~\ref{ass.AKP.2.1}~(a), and for some $\vep >0$, the inequality $\E[\Theta^{\alpha+ \vep}\,h_1(\Theta)] < \infty$ is fulfilled and holds the asymptotic relation
\beam \label{eq.AKP.2.2*}
\PP[\Theta> x]=o[\bF(x)]\,,
\eeam
then we find that
\beam \label{eq.AKP.2.13}
\PP[\Theta\,X>x] \sim \E[\Theta^{\alpha}\,h_1(\Theta)]\,\PP[X>x]\,.
\eeam
\epr
	
Under the convergence in \eqref{eq.AKP.2.3}, relation \eqref{eq.AKP.2.13} is equivalent to 
\beam \label{eq.AKP.2.14}
\lim_{\xto} x\,\PP\left[\dfrac {\Theta\,X}{U_F(x)} >1 \right]= \E[\Theta^{\alpha}\,h_1(\Theta)]\,.
\eeam

\bre \label{A*}
In the previous proposition, if we assume that $\E[\Theta^{\alpha + \vep}] < \infty$, then obviously relation \eqref{eq.AKP.2.11} holds, and further by \eqref{eq.AKP.2.7}, the condition $\E[\Theta^{\alpha + \vep}\,h_1(\Theta)] < \infty$ is satisfied. However, the conditions used in Proposition \ref{pr.AKP.2.1} are slightly weaker than  $\E[\Theta^{\alpha + \vep}] < \infty$, see \cite[Exam. 2.1]{cui:wang:2025} for discussions on the topic.  
\ere
	
The following assumption was introduced in \cite[Ass. 3.2]{konstantinides:leipus:passalidis:siaulys:2024} and relies on a bi-variate version of the dependence structure in Assumption~\ref{ass.AKP.2.1}. For the sake of brevity, we illustrate this asymptotic condition only for nonnegative random variables as in Assumption~\ref{ass.AKP.2.2}.
	
\begin{assumption} \label{ass.AKP.2.2}
Let $(X_i,\,Y_i)$ and $(\Theta_i,\,\Delta_i)$, $i \in \bbn$, be random vectors.  There exists a measurable function $g_i\;:\;[0,\,\infty)^2 \to (0,\,\infty)$, such that we have for any $i=1,\,\ldots,\,n$
\beam \label{eq.AKP.2.15}
\PP[X_i>x_1,\;Y_i>x_2\;\big|\;\Theta_i =\theta_i,\;\Delta_i=\delta_i] \sim g_{i}(\theta_i,\,\delta_i)\,\PP[X_i>x_1,\;Y_i>x_2]\,,
\eeam
as $x_1 \wedge x_2 \to \infty$, uniformly for any $(\theta_i,\;\delta_i) \in (S_{\Theta_i},\;S_{\Delta_i})$.
\end{assumption} 
	
Note that the uniform convergence in \eqref{eq.AKP.2.15} is understood as 
\beao
\lim_{x_1 \wedge x_2 \to \infty} \sup_{(\theta_i,\,\delta_i) \in (S_{\Theta},\,S_{\Delta})} \left| \dfrac{\PP[X_i > x_1,\;Y_i>x_2\;\big|\; \Theta_i=\theta_i,\;\Delta_i=\delta_i]}{g_{i}(\theta_i,\,\delta_i)\,\PP[X_i>x_1,\;Y_i>x_2]} -1 \right|=0\,,
\eeao
where, as with Assumption \ref{ass.AKP.2.1}, here again the left member of relation \eqref{eq.AKP.2.15} is understood in the sense
\beao
\lim_{\vep \downarrow 0} \PP\big[X_i > x_1\,,\;Y_i> x_2\;|\;\Theta_i \in [\theta_i - \vep,\,\theta_i + \vep]\,,\;\Delta_i \in [\delta_i - \vep,\,\delta_i + \vep]\big]\,.
\eeao 

\bre \label{rem.AKP.2.2}
It is easy to see that Assumption \ref{ass.AKP.2.2} contains the independence between the two random vectors $(X_i,\,Y_i)$ and $(\Theta_i,\,\Delta_i)$, as special case, while the pairs $(X_i,\,Y_i)$ have arbitrarily  dependent components, exactly as the pairs  $(\Theta_i,\,\Delta_i)$ as well. From \cite[Rem. 3.2]{konstantinides:leipus:passalidis:siaulys:2024} follows that 
\beao
\E[g_{i}(\Theta_i,\;\Delta_i)]=1\,,
\eeao
for any $i \in \bbn$ and further from \cite[Rem. 3.4]{konstantinides:leipus:passalidis:siaulys:2024}, there exists constant $\Lambda_i \in (0,\,\infty)$, such that the inequality
\beam \label{eq.KP.2.18}
g_{i}(\theta_i,\;\delta_i)\leq \Lambda_i\,,
\eeam
hold, for any $i \in \bbn$, for any $\theta_{i}\in S_{\Theta_{i}}$, $\delta_{i}\in S_{\Delta_{i}}$.
	\ere

\section{Breiman's theorem} \label{sec.KP.3}

	We now present an extension of Breiman's theorem in a complex dependent model that according to our knowledge, is the most general setting that exists in the literature; that is, the $BRV$ setting is assumed with dependent random weights that also dependent on the main random variables. One strand of research in the literature relies on the standard $MRV$ framework; e.g., see \cite{basrak:davis:mikosch:2002} and \cite{fougeres:mercadier:2012} for similar approaches  under either independence and dependence between random weights and their main random variables, respectively. Extensions beyond the non-standard $MRV$ are restricted to a two dimensional setting with independence between random weights and the main random variables; e.g., see \cite[Lem. 5.2]{chen:yang:2019}. Therefore, we present a generalization of this last reference that is the closest to our very general dependence framework. The following result, given as Theorem~\ref{th.AKP.3.1}, relies on 1) Assumption~\ref{ass.AKP.2.1} with \eqref{eq.AKP.2.9} and \eqref{eq.AKP.2.10} being satisfied with functions $h_1$ and $h_2$ for $(X,\,\Theta)$ and $(Y,\,\Delta)$ respectively, and 2) Assumption~\ref{ass.AKP.2.2} with \eqref{eq.AKP.2.15} being satisfied with function $g$ for $(X,\,\Theta)$ and $(Y,\,\Delta)$.
	
\bpr \label{th.AKP.3.1}
Let $(X,\,Y)$ and $(\Theta,\,\Delta)$ be nonnegative random vectors such that the pair $(X,\,Y) \in BRV_{-\alpha,\,-\beta}(\mu)$ for some finite $\alpha, \,\beta >0$, and Assumptions~\ref{ass.AKP.2.1} and \ref{ass.AKP.2.2} hold. Further, if we assume that 
\beao
\E[\Theta^{\alpha+\vep}\,h_1(\Theta)]+\E[\Delta^{\beta + \vep}\,h_2(\Delta)]< \infty\,,
\eeao 
for some $\vep >0$, and the relations
\beam \label{eq.AKP.2.18}
\PP[\Theta >x]=o[\bF(x)]\,, \qquad \PP[\Delta >x]=o[\bG(x)]\,,
\eeam
are valid, then, $(\Theta\,X,\,\Delta\,Y) \in BRV_{-\alpha,\,-\beta}(\widehat{\mu})$ with a limit measure $\widehat{\mu}$ satisfying for any $(p,\,q)> {\bf 0}$
\beam \label{eq.AKP.2.19} 
\widehat{\mu}\Big(\big([0,\,p)\times [0,\,q)\big)^c\Big) = \dfrac{\E[\Theta^{\alpha}\,h_1(\Theta)]}{p^\alpha}+\dfrac{\E[\Delta^{\beta}\,h_2(\Delta)]}{q^\beta} - \E\left[ g(\Theta,\,\Delta)\,\overline{\mu}\left( \dfrac p{\Theta},\;\dfrac q{\Delta} \right) \right]\,,
\eeam
where the tail of the two-dimensional measure $\overline{\mu}$ is given by
\beao
\overline{\mu}\left( x,\;y \right) :=\mu \left( \left(x,\,\infty \right]\times\left(y,\,\infty\right]\right) \,,\,\text{for any $x,\,y>0$.}
\eeao
\epr 
	
\pr~
For any $(p,\,q)> {\bf 0}$ by the assumptions on the moments and relation \eqref{eq.AKP.2.18}, we can apply Proposition \ref{pr.AKP.2.1}, and employing the formulation in relation \eqref{eq.AKP.2.14}, we obtain that
\beam \label{eq.AKP.2.20} 
\lim x\,\PP\left[\dfrac{\Theta\,X}{U_F(x)}>p \right]=\dfrac{\E[\Theta^{\alpha}\,h_1(\Theta)]}{p^\alpha}\,, \qquad \lim x\,\PP\left[\dfrac{\Delta\,Y}{U_G(x)}>q\right]=\dfrac{\E[\Delta^{\beta}\,h_2(\Delta)]}{q^\beta}\,,
\eeam
Let now some constant $C>0$, then we write
\beam \label{eq.AKP.2.21}  \notag
&&x\,\PP\left[\dfrac{\Theta\,X}{U_F(x)}>p,\;\dfrac{\Delta\,Y}{U_G(x)}>q \right]= x\,\PP\left[\dfrac{\Theta\,X}{U_F(x)}>p,\;\dfrac{\Delta\,Y}{U_G(x)}>q,\;\{\Theta \leq C\}\cap\{\Delta \leq C\} \right]\\[2mm] 
&& + x\PP\left[\dfrac{\Theta\,X}{U_F(x)}>p,\,\dfrac{\Delta\,Y}{U_G(x)}>q,\,\{\Theta > C\}\cup\{\Delta > C\} \right]=:R_1(x,\,p,\,q)+R_2(x,\,p,\,q),
\eeam
	
For the second term $R_2(x,\,p,\,q)$ in \eqref{eq.AKP.2.21}, we have the inequality
\beam \label{eq.AKP.2.22}  \notag
&& R_2(x,\,p,\,q) \\[2mm]  \notag
&&\leq x\,\PP\left[\dfrac{\Theta\,X}{U_F(x)}>p,\;\dfrac{\Delta\,Y}{U_G(x)}>q,\;\Theta > C \right]+ x\,\PP\left[\dfrac{\Theta\,X}{U_F(x)}>p,\;\dfrac{\Delta\,Y}{U_G(x)}>q,\;\Delta > C \right]\\[2mm]  \notag
&&\leq x\,\PP\left[\dfrac{\Theta\,X}{U_F(x)}>p,\; \Theta > C \right]+ x\,\PP\left[\dfrac{\Delta\,Y}{U_G(x)}>q,\;\Delta > C \right] \\[2mm] 
&&\sim \dfrac{\E[\Theta^{\alpha}\,h_1(\Theta)\,{\bf 1}_{\{\Theta>C\}}]}{p^\alpha}+\dfrac{\E[\Delta^{\beta}\,h_2(\Delta)\,{\bf 1}_{\{\Delta>C\}}]}{q^\beta}\,,
\eeam
where at the last step we used \eqref{eq.AKP.2.20}. Hence, taking into account the moment conditions, we apply the dominating convergence theorem on the last terms of \eqref{eq.AKP.2.22} to find that
\beam \label{eq.AKP.2.23}  
\lim_{\xto} R_2(x,\,p,\,q) = 0\,,
\eeam
as $C \to \infty$.
	
For the first term $R_1(x,\,p,\,q)$ in \eqref{eq.AKP.2.21}, we have 
\beam \label{eq.AKP.2.24}  \notag
&&\lim R_1(x,\,p,\,q) \\[2mm] \notag
&&= \lim \int_0^C \int_0^C x\,\PP\left[\dfrac{\theta\,X}{U_F(x)}>p,\;\dfrac{\delta\,Y}{U_G(x)}>q\;\Big|\;\Theta = \theta,\,\Delta = \delta \right]\,\PP[\Theta \in d\theta,\,\Delta \in d\delta]\\[2mm] \notag
&& = \int_0^C \int_0^C \lim x\,\PP\left[\dfrac{X}{U_F(x)}>\dfrac p{\theta},\;\dfrac{Y}{U_G(x)}>\dfrac q{\delta}\;\Big|\;\Theta = \theta,\,\Delta = \delta \right]\,\PP[\Theta \in d\theta,\,\Delta \in d\delta]\\[2mm] \notag
&&= \int_0^C \int_0^C \lim x\,g(\theta,\,\delta)\,\PP\left[\dfrac{X}{U_F(x)}>\dfrac p{\theta},\;\dfrac{Y}{U_G(x)}>\dfrac q{\delta} \right]\,\PP[\Theta \in d\theta,\,\Delta \in d\delta]\\[2mm] 
 &&= \int_0^C \int_0^C g(\theta,\,\delta)\,\overline{\mu}\left( \dfrac p{\theta},\;\dfrac q{\delta}\right) \,\PP[\Theta \in d\theta,\,\Delta \in d\delta] \to  \E \left[ g(\Theta,\,\Delta)\,\overline{\mu}\left( \dfrac p{\Theta},\;\dfrac q{\Delta}\right) \right]\,,
\eeam
as $C \to \infty$ at the last step, where at the third step we used Assumption \ref{ass.AKP.2.2} by the dominating convergence theorem. Thence, by relation \eqref{eq.AKP.2.21}, \eqref{eq.AKP.2.23} and \eqref{eq.AKP.2.24} we find that
\beam \label{eq.AKP.2.25} 
\lim_{\xto} x\,\PP\left[\dfrac{\Theta\,X}{U_F(x)}>p,\;\dfrac{\Delta\,Y}{U_G(x)}>q \right] = \E \left[ g(\Theta,\,\Delta)\,\overline{\mu}\left( \dfrac p{\Theta},\;\dfrac q{\Delta}\right) \right]\,.
\eeam
Finally, for any $(p,\,q)> {\bf 0}$ we obtain that
\beam \label{eq.AKP.2.26}
&&\lim_{\xto} x\,\PP\left[\left(\dfrac{\Theta\,X}{U_X(x)},\;\dfrac{\Delta\,Y}{U_Y(x)}\right) \in ([0,\,p)\times [0,\,q))^c \right] \\[2mm] \notag
&&= \lim_{\xto} x\,\Bigg(\PP\left[ \dfrac{\Theta\,X}{U_F(x)}>p\right] +\PP\left[\dfrac{\Delta\,Y}{U_G(x)}>q \right] -\PP\left[\dfrac{\Theta\,X}{U_F(x)}>p,\;\dfrac{\Delta\,Y}{U_G(x)}>q \right] \Bigg) \,.
\eeam
Now, combining the relations \eqref{eq.AKP.2.20}, \eqref{eq.AKP.2.25}, \eqref{eq.AKP.2.26}, we find \eqref{eq.AKP.2.19}. Furthermore, the $BRV$ structure for any Borel set $\bbb \subset \overline{\bbr}_{+}^2$, with ${\bf 0} \notin \overline{\bbb}$, follows by Dynkin's $\pi - \lambda$ theorem, see \cite[p. 36]{resnick:1999} or \cite[p. 178]{resnick:2007}. ~\halmos
	
\bre \label{rem.AKP.3.1}
In case where $(\Theta,\,\Delta)$ are independent of $(X,\,Y)$, in \cite[Lem. 5.2]{chen:yang:2019} was proved the previous proposition under the assumption that $\E\left[\Theta^{\alpha+\vep} \right] + \E\left[\Delta^{\beta+\vep} \right] < \infty$. Such an assumption, is implied by \eqref{eq.AKP.2.18} and furthermore from \eqref{eq.AKP.2.11} follows that $\E\left[\Theta^{\alpha+\vep}\,h_1(\Theta) \right] + \E\left[\Delta^{\beta+\vep}\,h_2(\Delta) \right] < \infty$. Hence, our result is immediately reduced to their result in the case of independence, and also provides weaker moment conditions in some concrete dependencies, as can be seen in \cite[Exam. 2.1]{cui:wang:2025}.
\ere

\section{BRV of finite randomly weighted sums} \label{sec.KP.4}

\subsection{Main results} \label{subsec.KP.4.1}

We now proceed with demonstrating the $BRV$ closure property with respect to randomly weighted sums $(S_n^{\Theta},\,T_n^{\Delta})$, for any fixed $n \in \bbn$. There exist many papers in the literature which study the distribution of the sum of heavy tailed random vectors; see \cite{meershaert:scheffler:2001}, \cite{hult:lindskog:mikosch:samorodnitsky:2005}, \cite{das:fashenhartmann:2023}, \cite{konstantinides:passalidis:2024b}, \cite{konstantinides:passalidis:2024g}, \cite{konstantinides:passalidis:2024h} among others.
	
In order to establish our main result, we need to present another assumption, that depicts the independence between the terms of the sequence  $\{ (X_i,\,Y_i,\,\Theta_i,\, \Delta_i)\,,\; i \in \bbn\}$. Assumption~\ref{ass.AKP.4.1} requires that the main random variables and the random weights do not depend for different time epochs. 
	
\begin{assumption} \label{ass.AKP.4.1}
$\{ (X_i,\,Y_i,\,\Theta_i,\,\Delta_i)\,,\; i \in \bbn\}$ is a sequence with independent terms, with respect to  $i \in \bbn$.
\end{assumption} 
	
At first glance, Assumption \ref{ass.AKP.4.1} could be restrictive in some models with insurance and financial risks, where $\Theta_i$ and $\Delta_i$ are 
usually understood as (random) discount factors over time and 
described as products of discount factors at each time period. However due to the fact that in Assumptions \ref{ass.AKP.2.1} and \ref{ass.AKP.2.2} the functions $h_{1,i},\,h_{2,i}$ and $g_i$, are NOT necessarily the same for all $i=1,\,\ldots,\,n$, the Assumption \ref{ass.AKP.4.1} finally seems not particularly restrictive for some credit risk applications, see in Subsec. $6.1$ for more details. Additionally, we adopt another assumption, related to the heavy-tailedness of the random weights  $\Theta_i$ and $\Delta_i$.
	
\begin{assumption} \label{ass.AKP.4.2}
For any $i \in \bbn$ we assume that 
\beao
\PP(\Theta_i >x ) =o[\bF(x)]\,, \qquad \PP(\Delta_i >x ) =o[\bG(x)]\,.
\eeao
\end{assumption} 
	
The following theorem is the main result fo this section.
	
\bth \label{th.AKP.4.1}
Let $\{(X_i,\,Y_i)\,,\;i \in \bbn\}$ be a sequence of i.i.d., nonnegative random vectors with a general random vector $(X,\,Y) \in BRV_{-\alpha,\,-\beta}(\mu)$ for some finite $\alpha,\,\beta >0$. Let $\{\Theta_i,\,\Delta_i\,,\;i \in \bbn\}$ be nonnegative and non-degenerate to zero, random variables. If Assumptions~\ref{ass.AKP.2.1}, \ref{ass.AKP.2.2}, \ref{ass.AKP.4.1} and \ref{ass.AKP.4.2} are satisfied and 
\beao
\E[\Theta_i^{\alpha+\vep}\,h_{1,i}(\Theta_i)]+\E[\Delta_i^{\beta + \vep}\,h_{2,i}(\Delta_i)]< \infty\,,
\eeao
for some $\vep >0$ and for all $i=1,\,\ldots,\,n$, then, 
\beao
(S_n^{\Theta},\,T_n^{\Delta}) \in BRV_{-\alpha,\,-\beta}(\widehat{\mu})\,,
\eeao 
for each $n \in \bbn$ with a limit measure $\widehat{\mu}$ satisfying 
\beam \label{eq.AKP.4.41}
&&\widehat{\mu}\left[\left([0,\,p) \times [0,\,q) \right)^c \right] \\[2mm]\notag
&&= \sum_{i=1}^n\Bigg(\dfrac{\E[\Theta_i^{\alpha}\,h_{1,i}(\Theta_i)]}{p^\alpha}+\dfrac{\E[\Delta_i^{\beta}\,h_{2,i}(\Delta_i)]}{q^\beta} - \E\left[ g_i(\Theta_i,\,\Delta_i)\,\overline{\mu}\left( \dfrac p{\Theta_i},\;\dfrac q{\Delta_i}\right) \right]\Bigg)\,,
\eeam 
for any $(p,\,q)> {\bf 0}$.
\ethe

We present two remarks that help in explanation of the degree of generality of the main results	of Theorem~\ref{th.AKP.4.1} (Remark~\ref{rem.AKP.4.3}), and to illustrate an important, direct application of  Theorem~\ref{th.AKP.4.1} ( Remark~\ref{rem.AKP.4.4}). 
	
\bre \label{rem.AKP.4.3}
The result in Theorem \ref{th.AKP.4.1}, looks similar with \cite[Th. 3.1]{chen:yang:2019}. However, there we find the assumption of the independence between the pairs $(X_i,\,Y_i)$ and $(\Theta_i,\,\Delta_i)$, for any $i \in \bbn$, but the $\{\Theta_i,\,\Delta_i\,,\; i \in \bbn\}$ are arbitrarily dependent. Here the pairs $(X_i,\,Y_i)$ and $(\Theta_i,\,\Delta_i)$ are weakly dependent (Assumptions \ref{ass.AKP.2.1}, \ref{ass.AKP.2.2}), but the $\{(\Theta_i,\,\Delta_i)\,,\; i\in \bbn \}$ are independent random vectors (with each pair $(\Theta_i,\,\Delta_i)$ possess arbitrary dependence). Hence, our result does not includes nor is included by that one. The applications that follow from our result (see Subsection $6.1$) are of different nature in comparison with the classical risk models found in the literature.
\ere
	
\bre \label{rem.AKP.4.4}	Investigating the behavior of concomitant extreme events is a fundamental application of Extreme Value Theory to risk management, but also its bespoke financial and insurance applications, see for example \cite{asimit:li:2018}, \cite{tang:xun:zhou:2026}.  Therefore, we examine the quantity
\beam \label{eq.AKP.4.56}
\PP[S_n^{\Theta}>p\,U_F(x)\;\big|\;T_n^{\Delta}>q\,U_G(x)]\,.
\eeam 
	
Under the conditions of Theorem~\ref{th.AKP.4.1}, one may find the following asymptotic formula
\beao
&&\lim \PP[S_n^{\Theta}>p\,U_F(x)\;\big|\;T_n^{\Delta}>q\,U_G(x)]\\[2mm] \notag
&&=\lim \dfrac {x\,\PP\left[\dfrac{S_n^{\Theta}}{U_F(x)}>p\,,\;\dfrac{T_n^{\Delta}}{U_G(x)}>q\right]}{x\,\PP\left[\dfrac{T_n^{\Delta}}{U_G(x)}>q\right]}=\dfrac {\displaystyle\sum_{i=1}^n\E\left[g_i(\Theta_i,\,\Delta_i)\,\overline{\mu}\left( \dfrac {p}{\Theta_i},\;\dfrac {q}{\Delta_i}\right) \right]}{\displaystyle\sum_{i=1}^n\E[\Delta_i^{\beta}\,h_{2,i}(\Delta_i)]}q^\beta\,,\notag
\eeao
where \eqref{eq.AKP.4.42} and \eqref{eq.AKP.4.55} are used at the last step.
\ere

\subsection{Argumentation} \label{subsec.KP.4.2}

Let us start with a preliminary lemma and some concepts, that are necessary for the proof of Theorem \ref{th.AKP.4.1}. We say that a distribution $V$, with $\bV(x)>0$ for any $x\in \bbr$, belongs to the class of dominatedly varying distributions, symbolically $V \in \mathcal{D}$, if for any (or, equivalently, for some) $b\in (0,\,1)$ holds
\beao
\limsup \dfrac {\bV(b\,x)}{\bV(x)} < \infty\,.
\eeao
It is easy to see that $\mathcal{R}:=\bigcup_{\alpha>0} \mathcal{R}_{-\alpha} \subsetneq \mathcal{D}$. Two important indexes for the distribution $V$ are the upper and lower Matuszewska indexes, that are defined for a distribution $V$, with $\bV(x)>0$ for any $x\in \bbr$, as follows
\beao
J^+_V:=\inf \left\{-\dfrac{\ln \bV_* (u)}{\ln u}\;:\;u>1 \right\}\,,\qquad J^-_V:=\sup \left\{-\dfrac{\ln \bV^* (u)}{\ln u}\;:\;u>1 \right\}\,,
\eeao
where
\beao
\bV_* (u):=\liminf \dfrac{\bV(ux)}{\bV(x)}\,,\qquad \bV^* (u):=\limsup \dfrac{\bV(ux)}{\bV(x)}\,.
\eeao
These indexes were introduced in \cite{matuszewska:1964} and satisfy the inequalities $0\leq J^-_V \leq J^+_V \leq \infty$. They characterize some well-known classes of heavy tailed distributions;  e.g., $V \in \mathcal{D}$ if and only if $J^+_V < \infty$, and further, if $V \in \mathcal{R}_{-\alpha}$  with $\alpha>0$ then $J^-_V=J^+_V=\alpha$; for details, see  \cite{bingham:goldie:teugels:1987} and \cite{leipus:siaulys:konstantinides:2023}.
	
The following lemma plays a crucial role in the proof of main result of this section given in Theorem~\ref{th.AKP.4.1}, but it also has its own merit as it represents a ``partial" dependent extension of  \cite[Lem. 7]{tang:yuan:2014}. Assumption ~\ref{ass.AKP.2.1}(a) below, holds in the sense of relation \eqref{eq.AKP.2.1*} for the pair $(X,\,\Theta)$. 
	
\ble \label{lem.AKP.4.1}
Let $X$ be a real-valued random variable with distribution $F \in \mathcal{D}$ and $\Theta$ be a nonnegative, non-degenerate to zero, random variable. Assume that $(X,\,\Theta)$ satisfies Assumption~\ref{ass.AKP.2.1}(a) and the moment condition $\E[\Theta^{\rho}\,h_1(\Theta)]< \infty$ for some $\rho> J^+_F$, and \eqref{eq.AKP.2.2*} holds as $\xto$ for any $c>0$. Further, let consider some set of random events $\{\Delta_t\,,\;t\in \mathcal{J} \}$, such that $\lim_{t\to t_0}\PP[\Delta_t] =0$ holds for some $t_0$ in the interior of set $\mathcal{J}$, with the $\{\Delta_t\,,\;t\in \mathcal{J} \}$ being independent of $X$. Then,
\beam \label{eq.AKP.2.27} 
\lim_{t\to t_0} \limsup_{\xto} \dfrac{\PP[\Theta\,X>x,\;\Delta_t]}{\PP[\Theta\,X>x]}=\lim_{t\to t_0} \limsup_{\xto} \dfrac{\PP[\Theta\,X>x,\;\Delta_t]}{\PP[X>x]}=0 \,.
\eeam
\ele
	
\pr~
Firstly, we show that the left fraction in relation \eqref{eq.AKP.2.27} converges to zero. If the random variable $\Theta$ has distribution $B$, define a new random variable $\Theta^*$ that is independent of $X$ and  $\{\Delta_t\,,\;t\in \mathcal{J} \}$ with a distribution given by
\beam \label{eq.AKP.2.28} 
B^*(d\theta)=h_1(\theta)\,B(d\theta) \,.
\eeam
Thus, $B^*$ becomes a proper distribution by keeping in mind \eqref{eq.AKP.2.12}. Further, since $\E[\Theta^{\rho}\,h_1(\Theta)]< \infty$ is true for some $\rho> J^+_F$, then for the same $\rho$ it follows that
\beam \label{eq.AKP.2.29}  
\E[(\Theta^*)^{\rho}] = \int_{0}^{\infty} \theta^{\rho}\,B^*(d\theta)=\int_{0}^{\infty} \theta^{\rho}\,h_1(\theta)\,B(d\theta) =\E[\Theta^{\rho}\,h_1(\Theta)]<\infty\,,
\eeam
where a constant $k_1>0$ exists due to the fact that function $h_1$ is bounded from above by keeping in mind \eqref{eq.AKP.2.11}. Inequality \eqref{eq.AKP.2.29}, the fact that $\Theta^*$ is independent of $X$ and since $X$ is independent of $\{\Delta_t\,,\;t\in \mathcal{J} \}$, we can apply \cite[Lem. 7]{tang:yuan:2014} to find that
\beam \label{eq.AKP.2.30} 
\lim_{t\to t_0} \limsup_{\xto} \dfrac{\PP[\Theta^*\,X>x,\;\Delta_t]}{\PP[\Theta^*\,X>x]}=0 \,.
\eeam
We shall show that the dependence can be solved in the products, namely it holds
\beam \label{eq.KP.c.1}
\PP(\Theta^*\,X > x) = \int_{0}^{\infty} h_1(s)\,\bF\left(\dfrac{x}{s} \right)\,B(ds) \sim \PP(\Theta\,X > x)\,.
\eeam
Initially, from \eqref{eq.AKP.2.29} and the inclusion $F \in \mathcal{D}$, through \cite[Th. 3.3 (iv)]{cline:samorodnitsky:1994} we find
\beam \label{eq.KP.c.2}
\PP(\Theta^*\,X > x) \asymp \bF\left(x \right)\,.
\eeam
From condition $F \in \mathcal{D}$ and the relation $\PP(\Theta >x ) =o[\bF(x)]$, we obtain from \cite[Prop. 3.1]{zhou:wang:wang:2012} that there exists positive function $a(\cdot)$, such that $a(x) \to \infty$, $a(x)= o(x)$ and 
\beam \label{eq.KP.c.3}
\PP(\Theta > a(x)) = o[\bF(x)]\,.
\eeam 
Consequently, 
\beao
&&\PP(\Theta\,X > x) = \left(\int_0^{a(x)} + \int_{a(x)}^{\infty}\right) \PP\left( X > \dfrac x{s} \;|\;\Theta =s \right)\,\PP(\Theta \in ds)\\[2mm]
&&=\int_0^{a(x)} \PP\left( X > \dfrac x{s} \;|\;\Theta =s \right)\,\PP(\Theta \in ds)+ O\left[\PP(\Theta > a(x)) \right]\\[2mm]
&&\sim \int_0^{a(x)} h_1(s)\bF\left( \dfrac x{s}\right)\,\PP(\Theta \in ds)+ o\left[\bF(x) \right]\\[2mm]
&&= \left(\int_0^{\infty} - \int_{a(x)}^{\infty} \right) h_1(s)\bF\left( \dfrac x{s}\right)\,\PP(\Theta \in ds)+ o\left[\PP(\Theta^*\,X > x) \right]\\[2mm]
&&= [1-o(1)]\,\PP(\Theta^*\,X > x)+ o\left[\PP(\Theta^*\,X > x) \right]\,,
\eeao 
where at the third step we used  Assumption \ref{ass.AKP.2.1} (a) in combination with \eqref{eq.AKP.2.30} \eqref{eq.KP.c.3}, while at the fourth step we took into consideration \eqref{eq.KP.c.2}.

Next, from relation \eqref{eq.AKP.2.30}, \eqref{eq.KP.c.1} we obtain
\beam \label{eq.AKP.2.32} 
\lim_{t\to t_0} \limsup_{\xto} \dfrac{\PP[\Theta^*\,X>x,\;\Delta_t]}{\PP[\Theta\,X>x]}=0 \,.
\eeam
	
Now, we define a new set of events with the same index set $\mathcal{J}$. Let $\{\Delta_t^*\,,\;t\in \mathcal{J} \}$ be the new set, which is by definition independent of all the sources of randomness, i.e., independent of $\Theta$, $\Theta^*$, $X$, $\{\Delta_t\,,\;t\in \mathcal{J} \}$, and this new set has the following properties
	\begin{enumerate}
	
	\item
	For some set $\mathcal{I}\subseteq\mathcal{J}$, with $t_0\in\mathcal{I}$, holds
	\beao
	\lim_{t\to z} \PP[\Delta_t]=\lim_{t\to z} \PP[\Delta_t^*]\,.
	\eeao
	for any $z\in\mathcal{I}$ and 
	\item
	\beao 
\lim_{t\to t_0} \limsup_{\xto} \PP[\Theta\,X>x,\;\Delta_t]\leq \lim_{t\to t_0} \limsup_{\xto} \PP[\Theta\,X>x,\;\Delta_t^*]\,.
	\eeao
	\end{enumerate}
	Then, 
	\beam \label{eq.AKP.2.33} \notag
\PP[\Theta^*\,X>x,\;\Delta_t] &=& \int_{0}^{\infty} h_1(\theta)\,\bF\left(\dfrac x{\theta} \right)\,B(d\theta)\,\PP[\Delta_t] \\[2mm] \notag
&\sim& \PP[\Delta_t^*]\,\int_{0}^{\infty} h_1(\theta)\,\bF\left(\dfrac x{\theta} \right)\,B(d\theta) \sim \PP[\Delta_t^*]\,\PP[\Theta\,X > x]\\[2mm]
&=&\PP[\Theta\,X>x,\;\Delta_t^*] \gtrsim \PP[\Theta\,X>x,\;\Delta_t]\,,
	\eeam
	with the asymptotic relations to hold as $t\to t_0$ and $\xto$, where at the first step we used \eqref{eq.AKP.2.28} and the fact that $\Theta^*$ independent of  $X$ and $\{\Delta_t\,,\;t\in \mathcal{J} \}$, at the second step the property (1) of  $\{\Delta_t^*\,,\;t\in \mathcal{J} \}$, at the third step the \eqref{eq.KP.c.1}, and at the last two steps the fact that $\{\Delta_t^*\,,\;t\in \mathcal{J} \}$, is independent of all the sources of randomness, and the property (2). Hence, from \eqref{eq.AKP.2.32} and \eqref{eq.AKP.2.33} we obtain the desired convergence
\beam \label{eq.AKP.2.34} 
\lim_{t\to t_0} \limsup_{\xto} \dfrac{\PP[\Theta\,X>x,\;\Delta_t]}{\PP[\Theta\,X>x]}=0 \,.
\eeam  
	
Secondly, we show that the second fraction in relation \eqref{eq.AKP.2.27} tends to zero. Because of \eqref{eq.KP.c.2} and together with \eqref{eq.KP.c.1} imply
\beam \label{eq.AKP.2.35} 
\PP[\Theta\,X>x]\asymp \PP[X>x] \,,
\eeam
Then \eqref{eq.AKP.2.34} and \eqref{eq.AKP.2.35} imply
\beao
\lim_{t\to t_0} \limsup_{\xto} \dfrac{\PP[\Theta\,X>x,\;\Delta_t]}{\PP[X>x]}=0 \,.~\halmos
\eeao
	
\bre \label{rem.AKP.4.1}
We should mention, that in case
\beao
F \in \mathcal{R}_{-\alpha}\,,
\eeao 
for some $\alpha>0$ (hence $J^+_F=\alpha$) and we use relation \eqref{eq.AKP.2.27} in  Theorem \ref{th.AKP.4.1} below as follows: the set of events $\{\Delta_x \}$ has index set  $\bbr_+$, that means we have the set $\{\Delta_x\,,\;x\in \bbr_+ \}$, that satisfies the convergence 
\beao
\lim_{\xto} \PP[\Delta_x]=0\,,
\eeao 
and relation \eqref{eq.AKP.2.27} takes the form
\beam \label{eq.AKP.2.36} 
\lim_{\xto} x\,\PP\left[\dfrac{\Theta\,X}{U_F(x)}>1,\;\Delta_x\right] = 0\,,
\eeam
see in \cite[Lem. 5.1]{chen:yang:2019}, for similar to \eqref{eq.AKP.2.36} approaches in independent case.
\ere
	
{\bf Proof of Theorem \ref{th.AKP.4.1}}~
Let $(p,\,q) > {\bf 0}$. From Assumptions~\ref{ass.AKP.2.1}, \ref{ass.AKP.4.1} and \ref{ass.AKP.4.2}, using \cite[Th. 1]{yang:wang:leipus:siaulys:2013} and Proposition \ref{th.AKP.3.1}, we obtain 
\beam \label{eq.AKP.4.42} \notag
&&\lim_{\xto} x\,\PP\left[\dfrac{S_n^{\Theta}}{U_F(x)}>p\right]= \sum_{i=1}^n\dfrac{\E[\Theta_i^{\alpha}\,h_{1,i}(\Theta_i)]}{p^\alpha}\,, \\[2mm]
&&	\lim_{\xto} x\,\PP\left[\dfrac{T_n^{\Delta}}{U_G(x)}>q\right]= \sum_{i=1}^n\dfrac{\E[\Delta_i^{\beta}\,h_{2,i}(\Delta_i)]}{q^\beta}\,,
\eeam
where we used the equivalent formulation of regular variation (see \eqref{eq.AKP.2.14}). For some fixed $\vep \in (0,\,1)$ we define the events
\beao
&&H_1(x):=\left\{ \bigvee_{i=1}^n \dfrac{\Theta_i\,X_i}{U_F(x)}> (1-\vep)\,p \right\}\,, \\[2mm] 
&& H_2(x):=\left\{ \bigvee_{i=1}^n \dfrac{\Delta_i\,Y_i}{U_G(x)}> (1-\vep)\,q \right\}\,,
\eeao 
and consequently, we find that 
\beam \label{eq.AKP.4.43} \notag
&&x\,\PP\left[\dfrac{S_n^{\Theta}}{U_F(x)}>p\,,\;\dfrac{T_n^{\Delta}}{U_G(x)}>q\right]= x\,\PP\left[\dfrac{S_n^{\Theta}}{U_F(x)}>p\,,\;\dfrac{T_n^{\Delta}}{U_G(x)}>q\,,\;H_1(x)\cap H_2(x)\right]\\[2mm] \notag
&& + x\,\PP\left[\dfrac{S_n^{\Theta}}{U_F(x)}>p\,,\;\dfrac{T_n^{\Delta}}{U_G(x)}>q\,,\;[H_1(x)\cap H_2(x)]^c\right]\\[2mm]
&&\leq x\,\PP\left[\dfrac{S_n^{\Theta}}{U_F(x)}>p\,,\;\dfrac{T_n^{\Delta}}{U_G(x)}>q\,,\;H_1(x)\cap H_2(x) \right] \\[2mm] \notag
&&+ x\,\PP\left[\dfrac{S_n^{\Theta}}{U_F(x)}>p\,,\;\dfrac{T_n^{\Delta}}{U_G(x)}>q\,,\;H_1^c(x)\right] +x\,\PP\left[\dfrac{S_n^{\Theta}}{U_F(x)}>p\,,\;\dfrac{T_n^{\Delta}}{U_G(x)}>q\,,\;H_2^c(x)\right]\\[2mm]\notag
&&=\sum_{k=1}^3 J_k(p,\,q,\,x) \,.
\eeam
Now, the first summand $J_1(p,\,q,\,x)$, satisfies
\beam \label{eq.AKP.4.44}\notag
&&J_1(p,\,q,\,x)\leq x\,\PP\left[H_1(x)\cap H_2(x)\right] \\[2mm] \notag
&&= x\,\PP\Bigg[\bigvee_{i=1}^n \dfrac{\Theta_i\,X_i}{U_F(x)}> (1-\vep)\,p\,,\;\bigvee_{j=1}^n \dfrac{\Delta_j\,Y_j}{U_G(x)}> (1-\vep)\,q \Bigg]\\[2mm] 
&&\leq \sum_{i=1}^n \sum_{j=1}^n x\,\PP\bigg[ \dfrac{\Theta_i\,X_i}{U_F(x)}> (1-\vep)\,p\,,\; \dfrac{\Delta_j\,Y_j}{U_G(x)}> (1-\vep)\,q \bigg] \\[2mm] \notag
&&=\left(\sum_{i=j=1}^n + \sum_{i\neq j=1}^n \right)x\,\PP\bigg[\dfrac{\Theta_i\,X_i}{U_F(x)}> (1-\vep)\,p\,,\; \dfrac{\Delta_j\,Y_j}{U_G(x)}> (1-\vep)\,q \bigg] =\sum_{l=1}^2 J_{1,l}(p,\,q,\,x) \,.
\eeam
For the first term $J_{1,1}(p,\,q,\,x)$, via \eqref{eq.AKP.2.25} is implied that
\beam \label{eq.AKP.4.47'} \notag
&&\lim_{\xto} x\,\PP\left[\dfrac{\Theta_i\,X_i}{U_F(x)}> (1-\vep)\,p\,,\; \dfrac{\Delta_i\,Y_i}{U_G(x)}> (1-\vep)\,q \right] \\[2mm] 
&&=\E\left[g_i(\Theta_i,\,\Delta_i)\,\overline{\mu}\left( \dfrac {(1-\vep)\,p}{\Theta_i},\;\dfrac {(1-\vep)\,q}{\Delta_i}\right) \right] \,.
\eeam
For the second term, $J_{1,2}(p,\,q,\,x)$, through the independence between $X_i$ and $\Delta_j\,Y_j$ for any $i\neq j$, by Assumption~\ref{ass.AKP.4.1} and due to moment condition, we can apply Lemma~\ref{rem.AKP.4.1} (see \eqref{eq.AKP.2.36}), from which we obtain
\beam \label{eq.AKP.4.46} \notag
\lim_{\xto} x\,\PP\left[\dfrac{\Theta_i\,X_i}{U_F(x)}> (1-\vep)\,p\,,\; \dfrac{\Delta_j\,Y_j}{U_G(x)}> (1-\vep)\,q \right] =0 \,.\\
\eeam
Thus, relations \eqref{eq.AKP.4.44}, \eqref{eq.AKP.4.47'} and \eqref{eq.AKP.4.46} suggest that 
\beam \label{eq.AKP.4.47} \notag
&&\limsup_{\xto} J_1(p,\,q,\,x) \leq  \sum_{i=1}^n \E\bigg[g_i(\Theta_i,\,\Delta_i)\,\overline{\mu}\left( \dfrac {(1-\vep)\,p}{\Theta_i},\;\dfrac {(1-\vep)\,q}{\Delta_i}\right) \bigg]\\[2mm] 
&&\leq  \sum_{i=1}^n\E\bigg[g_i(\Theta_i,\,\Delta_i)\,\overline{\mu}\left( \dfrac {(1-\vep)^{(\alpha \vee \beta)/\alpha}\,p}{\Theta_i},\;\dfrac {(1-\vep)^{(\alpha \vee \beta)/\beta}\,q}{\Delta_i}\right) \bigg]
\\[2mm]\notag
&&= (1-\vep)^{-(\alpha \vee \beta)} \sum_{i=1}^n\E\left[g_i(\Theta_i,\,\Delta_i)\,\overline{\mu}\left( \dfrac {p}{\Theta_i},\;\dfrac {q}{\Delta_i}\right) \right] \to \sum_{i=1}^n\E\left[g_i(\Theta_i,\,\Delta_i)\,\overline{\mu}\left( \dfrac {p}{\Theta_i},\;\dfrac {q}{\Delta_i}\right) \right]\,,
\eeam
as $\vep \to 0$, where the fact that $\mu$ is increasing function with respect to its component is used at the second step, while the homogeneity of measure $\mu$ is recalled at the third step.
	
Further, the second summand $J_2(p,\,q,\,x)$, satisfies
\beam \label{eq.AKP.4.48} \notag
J_2(p,\,q,\,x)&=& x\,\PP\left[\dfrac{S_n^{\Theta}}{U_F(x)}>p\,,\;\dfrac{T_n^{\Delta}}{U_G(x)}>q\,,\;H_1^c(x)\right]\\[2mm] \notag
	&\leq& x\,\PP\left[\dfrac{S_n^{\Theta}}{U_F(x)}>p\,,\;\bigvee_{j=1}^n\dfrac{\Theta_j\,X_j}{U_F(x)}\leq (1-\vep)\,p\right]\\[2mm] \notag
	&=&x\,\PP\Bigg[\dfrac{S_n^{\Theta}}{U_F(x)}>p\,,\;\bigvee_{i=1}^n\dfrac{\Theta_i\,X_i}{U_F(x)}> \dfrac pn\,,\;\bigvee_{j=1}^n\dfrac{\Theta_j\,X_j}{U_F(x)}\leq (1-\vep)\,p\Bigg]\\[2mm] 
	&\leq& \sum_{i=1}^n x\,\PP\left[\dfrac{\Theta_i\,X_i}{U_F(x)}> \dfrac pn\,,\;\sum_{i\neq l=1}^n\dfrac{\Theta_l\,X_l}{U_F(x)}> \vep\,p\right] \to 0\,,
\eeam
where the independence of $X_i$ and $\Theta_l\,X_l$, and Lemma \ref{lem.AKP.4.1} are used to conclude the last step. 

Similarly, via symmetrical arguments employed in the derivation of \eqref{eq.AKP.4.48}, it could help to show that the third summand $J_3(p,\,q,\,x)$, satisfies
	\beam \label{eq.AKP.4.49}
	J_3(p,\,q,\,x)  \to 0\,,\;\text{ as $\xto$.} 
	\eeam
	Combining \eqref{eq.AKP.4.43} and \eqref{eq.AKP.4.47} - \eqref{eq.AKP.4.49}, we conclude the upper bound result
\beam \label{eq.AKP.4.50}
\limsup_{\xto} x\,\PP\left[\dfrac{S_n^{\Theta}}{U_F(x)}>p\,,\;\dfrac{T_n^{\Delta}}{U_G(x)}>q\right] \leq \sum_{i=1}^n \E\left[g_i(\Theta_i,\,\Delta_i)\,\overline{\mu}\left( \dfrac {p}{\Theta_i},\;\dfrac {q}{\Delta_i}\right) \right] \,.
\eeam
	
Now, we examine the lower bound, for which we consider the events
\beao
E_i:=\left\{\dfrac{\Theta_i\,X_i}{U_F(x)}> p\,,\; \dfrac{\Delta_i\,Y_i}{U_G(x)}> q\right\}\,,
\eeao
for $i=1,\,\ldots,\,n$. Since $\Theta_i\,X_i$ and $\Delta_i\,Y_i$ are nonnegative for any  $i=1,\,\ldots,\,n$, Bonferroni inequality yields
\beam \label{eq.AKP.4.51}\notag
x\,\PP\left[\dfrac{S_n^{\Theta}}{U_F(x)}>p\,,\;\dfrac{T_n^{\Delta}}{U_G(x)}>q\right] \geq  x\,\PP\left[\bigcup_{i=1}^n E_i\right] &\geq& \sum_{i=1}^n x\,\PP[E_i]-\sum_{1\leq i<j \leq n} x\,\PP[E_i\cap E_j]\\[2mm]
&=:&I_1(p,\,q,\,x)- I_2(p,\,q,\,x)\,.
\eeam
The first term, $I_1(p,\,q,\,x)$, is dealt with the help of   \eqref{eq.AKP.2.25}, and we obtain that
\beam \label{eq.AKP.4.52}
\lim_{\xto} I_1(p,\,q,\,x)= \sum_{i=1}^n\E\left[g_i(\Theta_i,\,\Delta_i)\,\overline{\mu}\left( \dfrac {p}{\Theta_i},\;\dfrac {q}{\Delta_i}\right) \right]\,,
\eeam
while the second term, $I_2(p,\,q,\,x)$, satisfies 
\beam \label{eq.AKP.4.53}
\lim I_2(p,\,q,\,x)=0\,.
\eeam
Indeed, due to Lemma \ref{lem.AKP.4.1}, is implied that 
\beao
&&x\,\PP\bigg[ \dfrac{\Theta_i\,X_i}{U_F(x)}> p,\, \dfrac{\Delta_i\,Y_i}{U_G(x)}>q,\,\dfrac{\Theta_j\,X_j}{U_F(x)}> p,\; \dfrac{\Delta_j\,Y_j}{U_G(x)}>q \bigg] \\[2mm] \notag
&&\leq x \PP\left[\dfrac{\Theta_i X_i}{U_F(x)}> p,\, \dfrac{\Delta_j\,Y_j}{U_G(x)}>q \right]\to 0\,.
\eeao 
From relations ~\eqref{eq.AKP.4.51}, \eqref{eq.AKP.4.52}, \eqref{eq.AKP.4.53} we conclude that
\beam \label{eq.AKP.4.54}
\liminf_{\xto} x\,\PP\left[\dfrac{S_n^{\Theta}}{U_F(x)}>p\,,\;\dfrac{T_n^{\Delta}}{U_G(x)}>q\right] \geq \sum_{i=1}^n \E\left[g_i(\Theta_i,\,\Delta_i)\,\overline{\mu}\left( \dfrac {p}{\Theta_i},\;\dfrac {q}{\Delta_i}\right) \right]\,,
\eeam
in turn, \eqref{eq.AKP.4.50} and \eqref{eq.AKP.4.54} lead to
\beam \label{eq.AKP.4.55}
\lim_{\xto} x\,\PP\left[\dfrac{S_n^{\Theta}}{U_F(x)}>p\,,\;\dfrac{T_n^{\Delta}}{U_G(x)}>q\right] = \sum_{i=1}^n \E\left[g_i(\Theta_i,\,\Delta_i)\,\overline{\mu}\left( \dfrac {p}{\Theta_i},\;\dfrac {q}{\Delta_i}\right) \right]\,.
\eeam
Finally, relations \eqref{eq.AKP.4.42} and \eqref{eq.AKP.4.55} and the fact that
\beao
&&\lim_{\xto} x\,\PP\left[\left(\dfrac{S_n^{\Theta}}{U_F(x)}\,,\;\dfrac{T_n^{\Delta}}{U_G(x)}\right) \in ([0,\,p)\times [0,\,q))^c\right]\\[2mm] \notag
&&= \lim_{\xto} x\,\Bigg(\PP\left[\dfrac{S_n^{\Theta}}{U_F(x)}>p\right]+\PP\left[\dfrac{T_n^{\Delta}}{U_G(x)}>q\right]-\PP\left[\dfrac{S_n^{\Theta}}{U_F(x)}>p\,,\;\dfrac{T_n^{\Delta}}{U_G(x)}>q\right]\Bigg)\,,
\eeao
conclude \eqref{eq.AKP.4.41}, which completes the proof.
~\halmos

\section{$BRV$ of infinite randomly weighted sums } \label{sec.KP.5}

In this section we extend the validity of Theorem \ref{th.AKP.4.1} uniformly with respect to $n \in \bbn$ under some harder moment conditions, as these used there.

\subsection{Main result} \label{subsec.KP.5.1}

The following theorem is the main result of this section. Recall that for any $p \in (0,\,\infty)$ from $\E[\Theta^p] < \infty$ is implies that $\E[\Theta^p\,h_1(\Theta)] < \infty$, see also Remark \ref{rem.AKP.3.1}.

\bth \label{th.KP.5.1}
Let $\{(X_i,\,Y_i)\,,\;i \in \bbn\}$ be a sequence of i.i.d. nonnegative random vectors with general random vector $(X,\,Y) \in BRV_{-\alpha,\,-\beta}(\mu)$, with $\alpha,\,\beta > 0$. Let $\{\Theta_i,\,\Delta_i\,,\;i \in \bbn\}$ be a sequence of nonnegative and non-degenerate to zero random variables. We assume that Assumptions \ref{ass.AKP.2.1}, \ref{ass.AKP.2.2}, \ref{ass.AKP.4.1} are true. Further, we suppose that for some $\vep>0$ it holds
\beam \label{eq.AKP.5.1}
\sum_{i=1}^{\infty} \E\left[\Theta_i^{\alpha - \vep} \bigvee \Theta_i^{\alpha + \vep}  \right] < \infty\,, \qquad \sum_{i=1}^{\infty} \E\left[\Delta_i^{\beta - \vep} \bigvee \Delta_i^{\beta + \vep}  \right] < \infty\,.
\eeam
Then for any $n \in \bbn \cup \{\infty\}$ it holds
\beao
(S_n^{\Theta}\,,\;T_n^{\Delta}) \in BRV_{-\alpha,\,-\beta}(\widehat{\mu})\,,
\eeao
with measure $\widehat{\mu}$ given uniformly for $n \in \bbn$ by the relation
\beam \label{eq.AKP.5.2}
&&\widehat{\mu}\left[ ({[0,\,p) \times [0,\,q)})^c \right] := \\[2mm] \notag
&&\sum_{i=1}^{n} \,\left( \dfrac{\E\left[\Theta_i^{\alpha }\,h_{1,i}(\Theta_i)  \right]}{p^{\alpha}} +
\dfrac{\E\left[\Delta_i^{\beta }\,h_{2,i}(\Delta_i)  \right]}{q^{\beta}} - 
 \E\left[g_i(\Theta_i,\,\Delta_i)\,\overline{\mu}\left( \dfrac {p}{\Theta_i},\;\dfrac {q}{\Delta_i} \right) \right] \right)\,,
\eeam 
for any $(p,\,q) > {\bf 0}$.
\ethe

\bre \label{rem.KP.5.1}
We notice that from the moment condition of \eqref{eq.AKP.5.1}, Assumption \ref{ass.AKP.4.2} is implied immediately. As in Section $4$, in Remark \ref{rem.AKP.4.1}, here also the result can be compared with \cite[Th. 3.2]{chen:yang:2019}. But again, our result does not cover nor is covered by that. However, the focus there was on the case $n=\infty$, while ours is on uniformity over $n \in \bbn$. The uniformity of Theorem \ref{th.KP.5.1} except its mathematical elegance, is important in actuarial practice, since provides 'good' approximations in problems over infinite time horizon, via the finite horizon. We refer the reader to \cite[Sec. 3]{tang:2004b}, for more discussions about the value of the uniformity in risk theory.
\ere

\bre \label{rem.KP.5.1,5}
From the results of Theorems \ref{th.AKP.4.1} and \ref{th.KP.5.1}, we can see immediately that it holds
\beam \label{eq.AKP.5.2,5}
\lim x\,\PP\left(\left[\dfrac{S_n^{\Theta}}{U_F(x)}\,,\;\dfrac{T_n^{\Delta}}{U_G(x)}\right] \in \bbb \right)=\sum_{i=1}^n \lim x\,\PP\left(\left[\dfrac{\Theta_i\,X_i}{U_F(x)}\,,\;\dfrac{\Delta_i\,Y_i}{U_G(x)}\right] \in \bbb \right)\,,
\eeam
for any $\bbb \subset \overline{\bbr}_{+}^2$, with ${\bf 0} \notin \overline{\bbb}$, which is $\mu$-continuous. The previous relation provides a non-standard version of the multivariate linear single big jump principle, that usually is represented by the relation
\beao
\PP\left(\left[S_n^{\Theta}\,,\;T_n^{\Delta}\right] \in x\,\bbb \right) \sim \sum_{i=1}^n \PP \left( \left[\Theta_i\,X_i\,,\;\Delta_i\,Y_i\right] \in x\,\bbb\right)\,,
\eeao
see \cite[Sec. 4]{konstantinides:passalidis:2024g} for more discussion. Degenerating the weights to unit, we can obtain a even more explicit expression, through \eqref{eq.AKP.5.2,5}, for this principle. 

From the other hand side, in \cite[Cor. 4.10]{samorodnitsky:sun:2016} was proved that when $(X,\,Y)$ follows multivariate subexponential distribution, symbolically $\mathcal{S}_{\mathscr{R}}$, this principle is valid, but the inverse is not true in general (see \cite[Rem. 4.11]{samorodnitsky:sun:2016}). Hence, an essential open question is if $BRV$ (in its non-standard form) belongs to $\mathcal{S}_{\mathscr{R}}$ (for the standard $BRV$ this is true by \cite[Prop. 4.14]{samorodnitsky:sun:2016}).
\ere

\subsection{Argumentation} \label{subsec.KP.5.2}

Before the proof of Theorem \ref{th.KP.5.1}, we provide two important preliminary lemmas.

\ble \label{lem.KP.5.1}
Let $\{X_i\,,\; i \in \bbn\}$ be a sequence of i.i.d. nonnegative random variables with common distribution $F \in \mathcal{R}_{-\alpha}$, with $\alpha>0$. We assume that the $\{\Theta_i\,,\; i \in \bbn\}$ is a sequence of independent nonnegative and non-degenerate to zero random variables. If the $\{(X_i,\,\Theta_i)\,,\;i \in \bbn\}$ satisfy the Assumption \ref{ass.AKP.2.1} (a), and for some $\vep >0$ it holds
\beam \label{eq.AKP.5.3}
\sum_{i=1}^{\infty} \,\E\left[\Theta_i^{\alpha - \vep}\bigvee \Theta_i^{\alpha + \vep}  \right]< \infty\,,
\eeam 
then we find that
\beam \label{eq.AKP.5.4}
\lim_{N \to \infty} \limsup \dfrac{\PP\left(\sum_{i=N+1}^{\infty} \Theta_i\,X_i > x \right)}{\bF(x)} =\lim_{N \to \infty} \limsup \dfrac{\sum_{i=N+1}^{\infty} \PP\left(\Theta_i\,X_i > x \right)}{\bF(x)} =0\,.
\eeam 
\ele 

\pr~
Let some fixed $d>1$ and some large enough $N_0 \in \bbn$, such that it satisfies
\beao
\sum_{i=N_0 +1}^{\infty} i^{-d} \leq 1\,.
\eeao
Hence, for any $N > N_0$ it holds
\beam \label{eq.AKP.5.5}
&& \PP\left(\sum_{i=N+1}^{\infty} \Theta_i\,X_i > x \right)\leq \PP\left(\sum_{i=N+1}^{\infty} \Theta_i\,X_i > x\,\sum_{i=N +1}^{\infty}  i^{-d} \right) \\[2mm] \notag
&&\leq  \PP\left(\bigcup_{i=N+1}^{\infty} \{\Theta_i\,X_i > x \,i^{-d} \}\right) \leq \sum_{i=N+1}^{\infty} \PP\left(\Theta_i\,X_i > x \,i^{-d}\right) \\[2mm] \notag
&&\leq \sum_{i=N+1}^{\infty} K_i \PP\left(\Theta_i\,X_i > x \right) \lesssim \check{K}\,\bF(x)\,\sum_{i=N+1}^{\infty}\,\E\left[\Theta_i^{\alpha}\,h_{1,i}(\Theta_i) \right] \,,
\eeam
where the $K_i \in (0,\,\infty)$ are derived by the (generalized) Potter inequalities, see in \cite[Lem. 1 (b)]{li:2018}, and 
\beao
\check{K} = \sup_{i\geq N+1} K_i\,.
\eeao 
At the last step we used Proposition \ref{th.AKP.3.1} through dominated convergence theorem, because of condition \eqref{eq.AKP.5.3}.

From \eqref{eq.AKP.5.3} and  \eqref{eq.AKP.5.5}, for any $\vep > 0$ there exists some $N_0' \in \bbn$, such that for any $N > N_0 \vee N_0'$ it holds
\beao
\PP\left(\sum_{i=N+1}^{\infty} \Theta_i\,X_i > x  \right) \lesssim \vep\, \check{K}\,\bF(x)\,.
\eeao
From the last relation, taking into consideration the arbitrary choice of $\vep>0$, we obtain that the first limit in \eqref{eq.AKP.5.4} tends to zero.

For the second limit we use a entirely similar (but easier) argument, to have the same outcome.
~\halmos

The next lemma extends the result of \eqref{eq.AKP.4.42} uniformly for $n \in \bbn$, and has its own merit.

\ble \label{lem.KP.5.2}
Under the conditions of Lemma \ref{lem.KP.5.1}, it holds
\beam \label{eq.AKP.5.6}
\lim \sup_{n \in \bbn} \left| \dfrac{\PP\left(\sum_{i=1}^{n} \Theta_i\,X_i > x \right)}{\bF(x)\,\sum_{i=1}^{n} \,\E\left[\Theta_i^{\alpha}\,h_{1,i}(\Theta_i)\right]} - 1 \right|=0\,.
\eeam 
\ele 

\pr~
From \eqref{eq.AKP.4.42}, for any $N \in \bbn$, it holds
\beam \label{eq.AKP.5.7}
\lim \sup_{n \leq N} \left| \dfrac{\PP\left(\sum_{i=1}^{n} \Theta_i\,X_i > x \right)}{\bF(x)\,\sum_{i=1}^{n} \,\E\left[\Theta_i^{\alpha}\,h_{1,i}(\Theta_i)\right]} - 1 \right|=0\,.
\eeam 
It remains to show \eqref{eq.AKP.5.6} for $n > N$. Let some fixed $d>1$. We can find some large enough $N_1$, such that it holds
\beao
q:= \sum_{i=N_1 +1}^{\infty} i^{-d} < 1\,.
\eeao
Hence for any $N > N_1$ we obtain
\beam \label{eq.AKP.5.8}
&&\PP\left(\sum_{i=1}^{\infty} \Theta_i\,X_i > x \right) \\[2mm] \notag
&&= \PP\left(\sum_{i=1}^{\infty} \Theta_i\,X_i > x\,q\,,\; \sum_{i=N+1}^{\infty} \Theta_i\,X_i \leq x\,q\right) +\PP\left(\sum_{i=1}^{\infty} \Theta_i\,X_i > x\,q\,,\; \sum_{i=N+1}^{\infty} \Theta_i\,X_i > x\,q\right) \\[2mm] \notag
&& \leq \PP\left(\sum_{i=1}^{N} \Theta_i\,X_i > x \right) + \PP\left(\sum_{i=N+1}^{\infty} \Theta_i\,X_i > x\,q \right)=: L_1(x,\,q,\,N) + L_2(x,\,q,\,N)\,.
\eeam 

For $L_2(x,\,q,\,N)$, for any $\vep > 0$ there exists $N_2 \in \bbn$, with $N > N_1\vee N_2$, such that for any $N > N_2$ it holds
\beam \label{eq.AKP.5.9}
&&L_2(x,\,q,\,N) = \PP\left(\sum_{i=N+1}^{\infty} \Theta_i\,X_i > x\,\sum_{i=N_1 +1}^{\infty} i^{-d}  \right)  \\[2mm] \notag
&&\leq \sum_{i=N+1}^{\infty} \PP\left( \Theta_i\,X_i > x\,i^{-d}\right)  \leq \check{K} \sum_{i=N+1}^{\infty} \,\PP\left( \Theta_i\,X_i > x \right) < \vep\,\bF(x)\,,
\eeam
where at the third step we used \cite[Lem. 1 (b)]{li:2018}, ($K_i \in (0,\,\infty)$) and 
\beao
\check{K}=\sup_{i\geq 1} K_i\,,
\eeao 
while at the last step applied relation \eqref{eq.AKP.5.4}.

For $L_1(x,\,q,\,N)$, we obtain
\beam \label{eq.AKP.5.10}
&&L_1(x,\,q,\,N) \sim \left( \sum_{i=1}^{N} \,\E\left[\Theta_i\,h_{1,i}(\Theta_i)\right] \right) \,\bF[(1-q)\,x] \\[2mm] \notag
&&\sim (1-q)^{-\alpha}\,\,\bF[x] \sum_{i=1}^{N} \,\E\left[\Theta_i\,h_{1,i}(\Theta_i)\right] \lesssim \bF(x)\,\sum_{i=1}^{\infty} \,\E\left[\Theta_i^{\alpha}\,h_{1,i}(\Theta_i)\right]\,,
\eeam
where at the first step we used \eqref{eq.AKP.4.42}, while at the last step we let $N$ tend to infinity, that implies $(1-q) \uparrow 1$. Hence from \eqref{eq.AKP.5.8} -  \eqref{eq.AKP.5.10} we find
\beam \label{eq.AKP.5.11}
&&\PP \left( \sum_{i=1}^{\infty} \Theta_i\,X_i > x\right) \lesssim \bF(x)\, \sum_{i=1}^{\infty} \,\E\left[\Theta_i^{\alpha}\,h_{1,i}(\Theta_i)\right]\,.
\eeam

Then for any $\vep_0 > 0$ we can find some $N_3 \in \bbn$, and for $N> \bigvee_{i=1}^3 N_i$, for all $n > N$ it holds
\beam \label{eq.AKP.5.12}
&&\PP \left( \sum_{i=1}^{n} \Theta_i\,X_i > x\right) \leq \PP \left( \sum_{i=1}^{\infty} \Theta_i\,X_i > x\right) \lesssim \bF(x)\, \sum_{i=1}^{\infty} \,\E\left[\Theta_i^{\alpha}\,h_{1,i}(\Theta_i)\right] \\[2mm] \notag
&& \leq \bF(x)\, \left( \sum_{i=1}^{n} + \sum_{i=N+1}^{\infty} \right)\,\E\left[\Theta_i^{\alpha}\,h_{1,i}(\Theta_i)\right] \leq (1+\vep_0)\,\bF(x)\, \sum_{i=1}^{n}\,\E\left[\Theta_i^{\alpha}\,h_{1,i}(\Theta_i)\right]\,,
\eeam
where at the second step we used \eqref{eq.AKP.5.11}, and at the last step took into consideration \eqref{eq.AKP.5.3}.

From the other hand side, for any $\vep_0^* > 0$, there exists some $N_4 \in \bbn$, such that for 
\beao
N >\bigvee_{i=1}^4 N_i\,,
\eeao 
and for all $n > N$, with $N> N_4$, it holds
\beam \label{eq.AKP.5.13} \notag
&&\PP \left( \sum_{i=1}^{n} \Theta_i\,X_i > x\right) \geq \PP \left( \sum_{i=1}^{N} \Theta_i\,X_i > x\right) \sim \bF(x)\, \sum_{i=1}^{N} \,\E\left[\Theta_i^{\alpha}\,h_{1,i}(\Theta_i)\right] \\[2mm]
&& \geq \bF(x)\, \left( \sum_{i=1}^{n} - \sum_{i=N+1}^{\infty} \right)\,\E\left[\Theta_i^{\alpha}\,h_{1,i}(\Theta_i)\right]  \\[2mm] \notag 
&&\geq (1-\vep_0^*)\,\bF(x)\, \sum_{i=1}^{n}\,\E\left[\Theta_i^{\alpha}\,h_{1,i}(\Theta_i)\right]\,,
\eeam
where at the second step we used \eqref{eq.AKP.4.42} and at the last step took into account \eqref{eq.AKP.5.3}. From relations \eqref{eq.AKP.5.12} and \eqref{eq.AKP.5.13} we obtain
\beam \label{eq.AKP.5.14}
\lim \sup_{n>N} \left| \dfrac{\PP \left( \sum_{i=1}^{n} \Theta_i\,X_i > x\right) }{\bF(x)\, \sum_{i=1}^{n}\,\E\left[\Theta_i^{\alpha}\,h_{1,i}(\Theta_i)\right]}-1 \right|=0\,.
\eeam
From \eqref{eq.AKP.5.7} and \eqref{eq.AKP.5.14} we finally get  \eqref{eq.AKP.5.6}.
~\halmos

\bre \label{rem.KP.5.2}
We notice that the statements of Lemma \ref{lem.KP.5.1} and Lemma \ref{lem.KP.5.2} remain true even under  Assumption  \ref{ass.AKP.4.2} and the condition 
\beao
\sum_{i=1}^{\infty} \,\E\left[\Theta_i^{\alpha-\vep}\,h_{1,i}(\Theta_i) \bigvee \Theta_i^{\alpha+\vep}\,h_{1,i}(\Theta_i)\right] < \infty\,,
\eeao
for some $\vep > 0$, instead of \eqref{eq.AKP.5.3}. This can be seen following the lines of proofs of the lemmas above. However, we keep \eqref{eq.AKP.5.3} for sake of compactness with the proof of Theorem \ref{th.KP.5.1}. 
\ere

\noindent{\bf Proof of Theorem \ref{th.KP.5.1}.}~
Let 
\beao
(p,\,q)> {\bf 0}\,.
\eeao 
For any fixed $N \in \bbn$, the desired result holds uniformly for any $n \leq N$, by application of Theorem \ref{th.AKP.4.1}. It remains to show that it holds uniformly for any 
\beao
n > N\,.
\eeao
Initially, by Lemma \ref{lem.KP.5.2} it holds uniformly for $n \in \bbn$ (and not only for $n > N$) the limit relations:
\beam \label{eq.AKP.5.15} \notag
&&\lim x\,\PP \left( \dfrac{S_n^{\Theta}}{U_F(x)} > p \right) =\dfrac{ \sum_{i=1}^{n}\,\E\left[\Theta_i^{\alpha}\,h_{1,i}(\Theta_i)\right]}{p^{\alpha}}\,, \\[2mm] 
&&\lim x\,\PP \left( \dfrac{T_n^{\Delta}}{U_G(x)} > q \right) =\dfrac{ \sum_{i=1}^{n}\,\E\left[\Delta_i^{\beta}\,h_{2,i}(\Delta_i)\right]}{q^{\beta}}\,,
\eeam
It remains to show that it holds uniformly for $n> N$ the relation
\beam \label{eq.AKP.5.16} 
\lim x\,\PP \left( \dfrac{S_n^{\Theta}}{U_F(x)} > p \,,\; \dfrac{T_n^{\Delta}}{U_G(x)} > q \right) = \sum_{i=1}^{n}\,\E\left[g_i(\Theta_i,\,\Delta_i)\,\overline{\mu}\left(\dfrac p{\Theta_i},\,\dfrac q{\Delta_i} \right) \right]\,.
\eeam

At first we shall show that
\beam \label{eq.AKP.5.17} 
\sum_{i=1}^{\infty}\,\E\left[g_i(\Theta_i,\,\Delta_i)\,\overline{\mu}\left(\dfrac p{\Theta_i},\,\dfrac q{\Delta_i} \right) \right]< \infty\,.
\eeam
Indeed, from \eqref{eq.KP.2.18}, and the properties of homogeneity and increase of measure $\mu$ we obtain
\beao
&&\sum_{i=1}^{\infty}\,\E\left[g_i(\Theta_i,\,\Delta_i)\,\overline{\mu}\left(\dfrac p{\Theta_i},\,\dfrac q{\Delta_i} \right) \right] \\[2mm] 
&&=\sum_{i=1}^{\infty} \int_0^{\infty}\,\int_0^{\infty}\, g_i(\Theta_i,\,\Delta_i)\,\overline{\mu}\left(\dfrac p{\theta},\,\dfrac q{\delta} \right) \,\PP\left(\Theta_i \in d\theta\,,\;\Delta_i \in d\delta \right) \\[2mm] 
&&\leq \sum_{i=1}^{\infty} \,\Lambda_i\,\E\left[\overline{\mu}\left(\dfrac p{\Theta},\,\dfrac q{\Delta} \right) \right] \leq \check{\Lambda} \sum_{i=1}^{\infty}\,\E\left[\mu\left( {( p/{\Theta},\,\infty ]} \times {(0,\,\infty]} \right) \right]  \\[2mm] 
&&=  \check{\Lambda}\, \mu\left(( 1,\,\infty ] \times (0,\,\infty] \right)\,p^{-\alpha}\,\sum_{i=1}^{\infty}\,\E\left[ \Theta_i^{\alpha}\right] < \infty\,, 
\eeao 
where 
\beao
\check{\Lambda}= \sup_{i\geq 1} \Lambda_i\,.
\eeao 
Let some $\vep \in (0,\,1)$. From the \eqref{eq.AKP.5.17} we can find some $N_1 \in \bbn$, with 
\beao
N >N_1\,,
\eeao 
such that for all $n>N$ it holds
\beam \label{eq.AKP.5.18} \notag
&&\lim x\,\PP \left( \dfrac{S_n^{\Theta}}{U_F(x)} > p \,,\; \dfrac{T_n^{\Delta}}{U_G(x)} > q \right) \geq \lim x\,\PP \left( \dfrac{S_N^{\Theta}}{U_F(x)} > p \,,\; \dfrac{T_N^{\Delta}}{U_G(x)} > q \right) \\[2mm] \notag
&&+\lim x\,\PP \left( \dfrac{S_n^{\Theta}}{U_F(x)} > p \,,\; \dfrac{T_n^{\Delta}}{U_G(x)} > q \right)  \\[2mm] \notag 
&&\geq \lim x\,\PP \left( \dfrac{S_N^{\Theta}}{U_F(x)} > p \,,\; \dfrac{T_N^{\Delta}}{U_G(x)} > q \right) \\[2mm] \notag 
&&= \sum_{i=1}^N \E\left[g_i(\Theta_i,\,\Delta_i)\,\overline{\mu}\left(\dfrac p{\Theta_i},\,\dfrac q{\Delta_i} \right) \right] \\[2mm] \notag 
&& \geq \left(\sum_{i=1}^n - \sum_{i=N+1}^{\infty} \right) \E\left[g_i(\Theta_i,\,\Delta_i)\,\overline{\mu}\left(\dfrac p{\Theta_i},\,\dfrac q{\Delta_i} \right) \right] \\[2mm] 
&& \geq (1-\vep)\,\sum_{i=1}^n  \E\left[g_i(\Theta_i,\,\Delta_i)\,\overline{\mu}\left(\dfrac p{\Theta_i},\,\dfrac q{\Delta_i} \right) \right]\,,
\eeam
where at the second step we used relation \eqref{eq.AKP.4.55}. By the arbitrary choice of $\vep>0$ is implied the lower bound for \eqref{eq.AKP.5.16}, for all $n>N$.

For the upper bound of \eqref{eq.AKP.5.16}, we consider some $\vep \in (0,\,1)$. Then for any $n>N$ it holds
\beam \label{eq.AKP.5.19} \notag
&&\lim x\,\PP \left( \dfrac{S_n^{\Theta}}{U_F(x)} > p \,,\; \dfrac{T_n^{\Delta}}{U_G(x)} > q \right) \\[2mm] \notag 
&& = \lim x\,\PP \left( \dfrac{S_n^{\Theta}}{U_F(x)} > p \,,\; \dfrac{T_n^{\Delta}}{U_G(x)} > q\,,\;  \dfrac{\sum_{i=N+1}^n \Theta_i\,X_i}{U_F(x)} \leq \vep\,p \,,\; \dfrac{\sum_{i=N+1}^n \Delta_i\,Y_i}{U_G(x)} \leq \vep\,q \right)  \\[2mm] 
&& + \lim x\,\PP \Bigg( \dfrac{S_n^{\Theta}}{U_F(x)} > p \,,\; \dfrac{T_n^{\Delta}}{U_G(x)} > q\,,\left\{  \dfrac{\sum_{i=N+1}^n \Theta_i\,X_i}{U_F(x)} > \vep\,p \right\}\;\bigcup\; \\[2mm] \notag 
&& \left\{\dfrac{\sum_{i=N+1}^n \Delta_i\,Y_i}{U_G(x)} \leq \vep\,q \right\} \Bigg)\leq \lim x\,\PP \left( \dfrac{S_N^{\Theta}}{U_F(x)} > (1- \vep)\,p \,,\; \dfrac{T_N^{\Delta}}{U_G(x)} > (1- \vep)\,q \right) \\[2mm] \notag 
&& + \lim x\,\PP \Bigg( \dfrac{\sum_{i=N+1}^{\infty} \Theta_i\,X_i}{U_F(x)} > \vep\,p \Bigg) + \lim x\,\PP \Bigg( \dfrac{\sum_{i=N+1}^{\infty} \Delta_i\,Y_i}{U_G(x)} > \vep\,q \Bigg)  =: \sum_{i=1}^3 K_i(p,q,x)\,.
\eeam
For $K_1(p,q,x)$, following similar way as in \eqref{eq.AKP.4.47} we obtain for any $n> N$
\beam \label{eq.AKP.5.20}  \notag 
&&K_1(p,q,x) =\sum_{i=1}^{N}  \E\left[g_i(\Theta_i,\,\Delta_i)\,\overline{\mu}\left(\dfrac {(1-\vep)\,p}{\Theta_i},\,\dfrac {(1-\vep)\,q}{\Delta_i} \right) \right]\\[2mm] \notag 
&&\leq \sum_{i=1}^{n}  \E\left[g_i(\Theta_i,\,\Delta_i)\,\overline{\mu}\left(\dfrac {(1-\vep)^{(\alpha \vee \beta)/\alpha}\,p}{\Theta_i},\,\dfrac {(1-\vep)^{(\alpha \vee \beta)/\beta}\,q}{\Delta_i} \right) \right] \\[2mm]
&&= (1-\vep)^{-(\alpha \vee \beta)} \,\sum_{i=1}^{n}  \E\left[g_i(\Theta_i,\,\Delta_i)\,\overline{\mu}\left(\dfrac {p}{\Theta_i},\,\dfrac {q}{\Delta_i} \right) \right] \,,
\eeam
where at the first step we used \eqref{eq.AKP.4.55}.

For  $K_2(p,q,x)$ by Lemma \ref{lem.KP.5.1} we find
\beam \label{eq.AKP.5.21}  
\lim_{N \to \infty} \limsup (\vep\,p)^{\alpha}\,x\,\PP\left( \dfrac{\sum_{i=N+1}^{\infty} \Theta_i\,X_i}{U_F(x)} > \vep\,p \right) =0\,,
\eeam 
and in similar way
\beam \label{eq.AKP.5.22}
\lim_{N \to \infty} \limsup (\vep\,q)^{\beta}\,x\,\PP\left( \dfrac{\sum_{i=N+1}^{\infty} \Delta_i\,Y_i}{U_G(x)} > \vep\,q \right) =0\,.
\eeam 

From \eqref{eq.AKP.5.19} - \eqref{eq.AKP.5.22} we obtain the upper bound of \eqref{eq.AKP.5.16} holds uniformly for $n> N$. This in combination with \eqref{eq.AKP.5.15} provides the desired uniformity for any $n > N$.
~\halmos

\section{Applications of our main results} \label{subsec.KP.6}

\subsection{Ruin probability in a credit risk model} \label{subsec.KP.6.1}

Now, we provide an application in a credit risk model and concretely we estimate the asymptotic behavior of three types of ruin probabilities. In the discrete-time risk model, the $(\Theta_i,\,\Delta_i)$ usually represent the stochastic discount coefficients, that makes the independence condition in the terms of the sequence 
\beao
\{ (\Theta_i,\,\Delta_i)\,,\; i\in \bbn \}\,, 
\eeao
to be hard (see, Assumption  \ref{ass.AKP.4.1}). Therefore, in the usual discrete time risk models, very often is used direction (i), or its similar extensions. However, in credit risk, the sequence $\{ (\Theta_i,\,\Delta_i)\,,\; i\in \bbn \}$ can represent the breach or non-breach of the obliged, while the sequence 
\beao
\{ (X_i,\,Y_i)\,,\; i\in \bbn \}\,,
\eeao 
may represent the amount of breach. To provide more flexibility to our model, we consider that the $\{ (X_i,\,Y_i)\,,\; i\in \bbn \}$ are real-valued (for example, if the realizations $\Theta_1 = 1$, and $X_1$ is negative, eventually means that the obliged paid after some delay his debt, and for this reason was billed more). Further, the $n$-obliged have a concrete risk profile (reflected by the $X_1,\,X_2,\,\ldots,\,X_n$), while the others $n$-obliged have different one (reflected by the $Y_1,\,Y_2,\,\ldots,\,Y_n$). We here have 
\beao
(X,\,Y) \in BRV_{-\alpha}(\mu)\,,
\eeao
namely this pair belongs to the standard $BRV$, as follows from the relation \eqref{eq.AKP.2.7} and that the function $U_V$ represents the normalizing function, for $V \in \mathcal{R}_{-\alpha}$, with finite $\alpha >0$.
	
The current results in the literature, related with ruin probabilities in discrete time models, focus mostly in the so-called insurance and financial risk models, where  Assumption \ref{ass.AKP.4.1}, used here, is heavy enough for the discount coefficients.

The most of that models focus in direction (i), of Section $1$, although some recent papers succeed to step in a common study of both directions (i) and (ii), see for example \cite{cheng:konstantinides:wang:2024}, but under the condition of asymptotic dependence to the corresponding $(\Theta\,X\,,\;\Delta\,Y)$, namely for 
\beam \label{eq.AKP.6.0}
\widehat{\mu}[(1,\,\infty]\times (1,\,\infty]]>0\,.
\eeam
Here, we do not make use of such condition, and thus we permit cases of asymptotic independence for $(X,\,Y)$. 

	Before depicting the risk model of this section, we provide a preliminary lemma, which is helpful in extension of our result in such pair $(X,\,Y)$ of random variables, that can take values over the whole real line.
	
	In what follows, the sums $S_n^{\Theta}, \,T_n^{\Delta}$, were introduced in relations \eqref{eq.AKP.1.1}, but now with $X_i,\,Y_i$ to be real-valued random variables, and we write
\beam \label{eq.AKP.6.1} 
S_n^{\Theta+}:=\sum_{i=1}^n \Theta_i\,X^{+}_i\,,\qquad \qquad T_n^{\Delta+}:=\sum_{i=1}^n \Delta_i\,Y^{+}_i\,.
\eeam
In what follows, when we say that the assumptions of Theorem \ref{th.AKP.4.1} holds, we mean that it holds for the positive part only, namely for the vector $(X^{+},\,Y^{+})$.
	
\ble \label{lem.AKP.5.1}
Under the conditions of Theorem \ref{th.AKP.4.1}, with the restriction that 
\beao
(X^{+},\,Y^{+})\in BRV_{-\alpha}(\mu)\,,
\eeao 
for finite $\alpha>0$, for any $(p,\,q) >{\bf 0}$, we obtain that
\beam \label{eq.AKP.6.2} \notag
\PP[S_n^{\Theta}> p\,x\,,\; T_n^{\Delta}>q\,x] &\sim&
\PP[S_n^{\Theta+}> p\,x\,,\; T_n^{\Delta+}>q\,x] \\[2mm] 
&\sim& \bV(x)\,\sum_{i=1}^n \E\left[g_i(\Theta_i,\,\Delta_i)\,\overline{\mu}\left( \dfrac {p}{\Theta_i},\;\dfrac {q}{\Delta_i}\right) \right]\,.
\eeam
\ele

\pr~
Initially from relation \eqref{eq.AKP.4.55} for any $p,\,q > 0$ and 
with the common normalizing function $U_V$, since we consider the standard $
BRV$, we find that
\beao
\lim x \PP\left[\dfrac{(S_n^{\Theta+}, T_n^{\Delta+})}{U_V(x)}\in (
p, \infty]\times (q, \infty]  \right] = \sum_{i=1}^n \E\left[g_i(\Theta_i, \Delta_i)\,\overline{\mu}\left( \dfrac {p}{\Theta_i},\;\dfrac {q}{\Delta_i}\right) \right
]\,,
\eeao
that after change of variables can be written as
\beam \label{eq.AKP.6.3}
\PP[S_n^{\Theta+}> p\,x\,,\; T_n^{\Delta+}>q\,x] \sim \bV(x)\,\sum_{i=1}^n \E\left[g_i(\Theta_i,\,\Delta_i)\,\overline{\mu}\left( \dfrac {p}{\Theta_i},\;\dfrac {q}{\Delta_i}\right) \right]\,.
\eeam
Obviously, from relations \eqref{eq.AKP.6.1} and \eqref{eq.AKP.6.3} we conclude that
\beam \label{eq.AKP.6.4} \notag
\PP[S_n^{\Theta}> p\,x\,,\; T_n^{\Delta}>q\,x]&\leq& \PP[S_n^{\Theta+}> p\,x
\,,\; T_n^{\Delta+}>q\,x] \\[2mm]
&\sim& \bV(x)\,\sum_{i=1}^n \E\left[g_i(\Theta_i,\,\Delta_i)\,\overline{\mu}\left( \dfrac {p}{\Theta_i},\;\dfrac {q}{\Delta_i}\right) \right]\,.
\eeam
Now, we take the lower bound in relation \eqref{eq.AKP.6.2}. 
Let fix some $\vep \in (0,\,1)$, then by Bonferroni inequality we obtain that
\beam \label{eq.AKP.6.5} \notag
&&\PP[S_n^{\Theta}> p\,x\,,\; T_n^{\Delta}>q\,x] \\[2mm] \notag
&&\geq \PP\Bigg[S_n^{\Theta}> p\,x\,,\; T_n^{\Delta}>q\,x\,,\; \bigcup_{i=1}^{n}\Big\{\Theta_i\,X_i>(1+\vep)\,p\,x\,,\Delta_i\,Y_i>(1+\vep)\,q\,x \Big\}\Bigg] \\[2mm] \notag
&&\geq \sum_{i=1}^n \PP \big[S_n^{\Theta}> p\,x\,,\; T_n^{\Delta}>q\,x\,,\; \Theta_i\,X_i>(1+\vep)\,p\,x\,,\Delta_i\,Y_i>(1+\vep)\,q\,x \big]\\[2mm] \notag  
&& -  \sum_{1\leq i<j\leq n}\PP\big[S_n^{\Theta}> p\,x\,,\; T_n^
{\Delta}>q\,x\,,\; \Theta_i\,X_i>(1+\vep)\,p\,x\,,\;\Delta_i\,Y_i>(1+\vep)\,q\,x\,,\\[2mm] 
&&\qquad \qquad \; \Theta_j\,X_j >(1+\vep)\,p\,x\,,\;\Delta_j\,Y_j>(1+\vep
)\,q\,x \big]\\[2mm] \notag
&&= \sum_{i=1}^n \PP\big[\Theta_i\,X_i>(1+\vep)\,p\,x\,,\;\Delta_i\,Y_i>(1+ \vep)\,q\,x \big]\\[2mm] \notag
&&- \sum_{i=1}^n \PP\big[\Theta_i\,X_i>(1+\vep)\,p\,x\,, \;\Delta_i\,Y
_i>(1+\vep)\,q\,x\,, \;\{S_n^{\Theta}> p\,x\,,\; T_n^{\Delta}>q\,x\}^c \big]\\[2mm]\notag
&&-\sum_{1\leq i<j\leq n}\PP\big[S_n^{\Theta}> p\,x\,,\; T_n^{\Delta}>q\,x\,,\;\Theta_i\,X_i>(1+\vep)\,p\,x\,,\;\Delta_i\,Y_i>(1+\vep)\,q\,x\,, \\[2mm] \notag
&&\qquad \qquad \;\Theta_j\,X_j>(1+\vep)\,p\,x\,,\;\Delta_j\,Y_j>(1+
\vep)\,q\,x \big] \\[2mm] \notag
&&=: \sum_{i=1}^3 K_i(p,\,q,\,x) \,.
\eeam
For the $K_3(p,\,q,\,x)$, since the $X_i$ and $\Delta_j\,Y_j$ are independent for $i<j$, taking into consideration that $	V \in \mathcal{R}_{-\alpha} \subsetneq \mathcal{D}$, 	with $\alpha> 0$ and  in view of the moment conditions of  the 	Lemma \ref{lem.AKP.4.1}, we obtain 
\beam \label{eq.AKP.6.6} 
K_3(p,\,q,\,x)\leq \sum_{1\leq i<j\leq n} \PP\left[\Theta_i\,X_i>(1+
\vep)\,p\,x\,,\;\Delta_j\,Y_j>(1+\vep)\,q\,x \right]=o[\bV(x)] \,.
\eeam
For the $K_1(p,\,q,\,x)$, using the fact that 
\beao
\PP\left[\Theta\,X>p\,x\,,\;\Delta\,Y>q\,x \right]=\PP\left[\Theta\,X^{+}>p\,x\,,\;\Delta\,Y^{+}>q\,x \right]\,,
\eeao
for any $x>0$ and any $(p,\,q) >{\bf 0}$, and by relation \eqref{eq.AKP.6.3} with $n=1$ for the standard $BRV$, we find that
\beao
K_1(p,\,q,\,x)&=&\sum_{i=1}^n\PP\left[\Theta_i\,X_i^{+}>(1+\vep)\,p\,x\,,
\;\Delta_i\,Y_i^{+}>(1+\vep)\,q\,x \right] \\[2mm] 
&\sim &\bV(x)\,\sum_{i=1}^n \E\left[g_i(\Theta_i,\,\Delta_i)
\,\overline{\mu}\left( \dfrac {(1+\vep)\,p}{\Theta_i},\;\dfrac {(1+\vep)\,q}{\Delta_i}\right) \right] \\[2mm] 
	&=&\dfrac{\bV(x)}{ (1+\vep)^{\alpha}}\,\sum_{i=1}^n \E\left[g_i(\Theta_i,\,\Delta_
	i)\,\overline{\mu}\left( \dfrac {p}{\Theta_i},\;\dfrac {q}{\Delta_i}\right) \right]\,,
	\eeao
where at the last step we used the property of homogeneity of measure $\mu$. Therefore we conclude that
	\beam \label{eq.AKP.6.7}
	K_1(p,\,q,\,x) 	\sim \bV(x)\,\sum_{i=1}^n \E\left[g_i(\Theta_i,\,\Delta_i)\,\overline{\mu}\left( \dfrac {p}{\Theta_i},\;\dfrac {q}{\Delta_i}\right) \right]\,,
	\eeam
	as $\vep \downarrow 0$. Further, we study $K_2(p,\,q,\,x)$ as follows 
	\beam \label{eq.AKP.6.8}
	&&K_2(p,\,q,\,x) \leq \sum_{i=1}^n \PP\big[\Theta_i\,X_i^{+}>(1+\vep)\,p\,x\,,\;\Delta_i\,Y_i^{+}>(1+\vep)\,q\,x\,,\;S_n^{\Theta}\leq p\,x \big] \\[2mm] \notag
	&&+ \sum_{i=1}^n\PP\big[\Theta_i\,X_i^{+}>(1+\vep)\,p\,x\,,\;\Delta_i\,Y_i^{+}>(1+\vep)\,q\,x\,,\;T_n^{\Delta} \leq q\,x\big]=:\sum_{j=1}^2 K_{2j}(p,\,q,\,x)\,.\notag
	\eeam
The first term can be bound as follows
\beam \label{eq.AKP.6.9} \notag
&& 	K_{21}(p,\,q,\,x)\leq\sum_{i=1}^n \PP\bigg[\Theta_i\,X_i^{+}>(1+\vep)\,p\,x\,,\; \sum_{i\neq k=1}^n \Theta_k\,X_k \leq -\vep\,p\,x \bigg]\\[2mm] 
&&= o\left(\PP\left[X_i^{+} >(1+\vep)\,x \right] \right) =o[\bV(x)]\,,
\eeam
where at the second step we used Lemma \ref{lem.AKP.4.1}, since $X_i^{+}$ is independent of all $\Theta_k\,X_k$, for any $i\neq k$, and further we employ either the fact that $\PP\left[X_i^{+}>(1+\vep)\,x\right] \asymp \bV[(1+\vep)\,x]$ and then we complete the argument, in view of the inclusion $V \in \mathcal{R}_{-\alpha}$.

Via symmetrical arguments we obtain
\beam \label{eq.AKP.6.10}
K_{22}(p,\,q,\,x)=o[\bV(x)]\,.
\eeam
Hence, from relations \eqref{eq.AKP.6.8}, \eqref{eq.AKP.6.9} and \eqref{eq.AKP.6.10} we conclude that
\beam \label{eq.AKP.6.11}
K_{2}(p,\,q,\,x)=o[\bV(x)]\,.
\eeam
Therefore from \eqref{eq.AKP.6.5} - \eqref{eq.AKP.6.7} and  \eqref{eq.AKP.6.11} we have the asymptotic inequality
\beam \label{eq.AKP.6.12}
\PP[S_n^{\Theta}> p\,x\,,\; T_n^{\Delta}>q\,x] \gtrsim \bV(x)\,\sum_{i=1}^n \E\left[g_i(\Theta_i,\,\Delta_i)\,\overline{\mu}\left( \dfrac {p}{\Theta_i},\;\dfrac {q}{\Delta_i}\right) \right] \,.
\eeam
Thus finally, by relations \eqref{eq.AKP.6.4} and \eqref{eq.AKP.6.12} we find \eqref{eq.AKP.6.2}. 
~\halmos
	
We remind that this section highlights the asymptotic estimation of the finite-time ruin probability, in a bivariate, discrete-time, credit risk model. More and more papers focus their attention on bivariate risk model, as a first step to explore the multivariate risk model and further to explore the dependence modeling among several parallel business lines in insurance industry, something by no doubt significant in actuarial practice, see for example \cite{chen:li:cheng:2023}, \cite{sun:yuan:lu:2024}, \cite{wang:su:yang:2024}, \cite{yang:chen:yuen:2024} among others.
	
We keep the representation of the credit-risk model as described at the beginning of the section (namely, $X_i,\,Y_i,\,\; i= 1,\,\ldots,\,n$ are the quantity of breach / or of gain by the delayed  payment of the obliged, and $\Theta_i,\,\Delta_i,\,\; i= 1,\,\ldots,\,n$ depict the breach (or, partial breach) or non-breach of the obliged). So, if the initial capital of the financial entity is $x>0$ and it is split to two lines in the form $p\,x$ and $q\,x$, where $p,\,q >0$ and $p+q=1$, then the surplus processes to each business line up to moment $k=1,\,\ldots,\,n$ (namely the moment where the $k$-th pair $(X_k,\,Y_k)$ enters into the system) can be written in the form
\beam \label{eq.AKP.6.13}
&&U_1(p\,x,\,k)=p\,x - \sum_{i=1}^k \Theta_i\,X_i\,,\\[2mm] \notag
&& U_2(q\,x,\,k)=q\,x - \sum_{j=1}^k \Delta_j\,Y_j\,.
\eeam 
The depiction of \ref{eq.AKP.6.13}, provides the flexibility that the $n$ pairs of customers are not coming nor leaving necessarily together to or from the system. This kind of flexibility is not given by several recent credit risk models, as for example \cite{chen:yang:zhang:2024}, \cite{chen:tong:yang:2025}, although in that models we find other kinds of flexibilities that are missing from our application. Furthermore, an other flexibility of our model, is that often in practice the credibility of an obliged (namely the $\Theta_i,\,\Delta_i$) are positively dependent with the quantity of breach (namely the $X_i,\,Y_i$).

In the one-dimensional set up the ruin is defined as the event when the surplus becomes negative, and consequently we find the one-dimensional ruin times up to the time horizon $n$ as follows
	\beao
	\tau_1(p\,x)&=&\inf\left\{k\in\mathbb{N}:\;U_1(p\,x,\,k)<0\;|\; U_1(p\,x,\,0)=p\,x \right\}\,, \\[2mm]
	\tau_2(q\,x)&=&\inf\left\{k\in\mathbb{N}:\;U_2(q\,x,\,k)<0\;|\; U_2(q\,x,\,0)=q\,x \right\}\,.
	\eeao
	
	In bivariate case, the ruin probability can be defined in several different ways, see \cite[Sec. 1]{cheng:yu:2019}. Here, we use three types of the ruin probability at finite horizon, up to moment $n$. The first type is the probability that both portfolios to fall below zero, up to moment $n$, but not necessarily simultaneously. Hence, this ruin time is defined as follows
\beao	
\tau_{and}=\max\left\{\tau_1(p\,x),\,\tau_2(q\,x) \right\}\,,
\eeao
and consequently the ruin probability from relation \eqref{eq.AKP.6.13} is defined as
\beao
\psi_{and}(x,\,p,\,q;\,n)=\PP[\tau_{and}\leq n]\,,
\eeao
or, equivalently, as
\beam \label{eq.AKP.6.14}
\psi_{and}(x,\,p,\,q;\,n) =\PP\left[\max_{1\leq k \leq n} \sum_{i=1}^k \Theta_i\,X_i >p\,x\,,\;\max_{1\leq l \leq n} \sum_{j=1}^l \Delta_j\,Y_j >q\,x\right]\,,
\eeam
	
The second kind of ruin probability, corresponds to the case when both portfolios fall simultaneously below zero. Thus, this ruin time is given by
\beao
\tau_{sim}=\inf\left\{k \in \bbn \;:\;U_1(p\,x,\,k)\vee U_2(p\,x,\,k)<0 \right\}\,,
\eeao 
and the ruin probability takes the form $\psi_{sim}(x,\,p,\,q;\,n)=\PP[\tau_{sim}\leq n]$.
	
Finally, the third kind of ruin probability, appears when at least one of the two portfolios falls below zero. In this case the ruin time is defined as
\beao
\tau_{or}=\min \left\{\tau_1(p\,x)\,,\;\tau_2(q\,x)\right\}\,,
\eeao 
with corresponding ruin probability of the form 
\beao
\psi_{or}(x,\,p,\,q;\,n)=\PP[\tau_{or} \leq n]= \PP\left[ \left\{\max_{1\leq k \leq n} \sum_{i=1}^k  \Theta_i X_i >p\,x\right\}\bigcup \left\{\max_{1\leq l \leq n} \sum_{i=1}^l  \Delta_i Y_i >q\,x\right\} \right].
\eeao
	
Now we proceed with the following theorem
	
\bth \label{th.AKP.5.1}
Assume that the conditions in Lemma \ref{lem.AKP.5.1} hold with $p,\,q>0$ and $p+q=1$.
\begin{enumerate}
\item[(i)] For each $n \in \bbn$, we obtain that 
\beam \label{eq.AKP.6.15}
\psi_{and}(x,\,p,\,q;\,n) \sim \psi_{sim}(x,\,p,\,q;\,n) \sim \bV(x) \sum_{i=1}^n \E\left[g_i(\Theta_i,\,\Delta_i)\,\overline{\mu}\left( \dfrac {p}{\Theta_i},\;\dfrac {q}{\Delta_i}\right) \right]\,.
\eeam 
\item[(ii)] For each $n \in \bbn$, we obtain that 
\beam \label{eq.AKP.6.17}
&& \psi_{or}(x,\,p,\,q;\,n)\sim \\[2mm] \notag
&&\bV(x)\,\Bigg( \sum_{i=1}^n\dfrac{\E[\Theta_i^{\alpha}\,h_{1,i}(\Theta_i)]}{p^\alpha} + \sum_{i=1}^n\dfrac{\E[\Delta_i^{\alpha}\,h_{2,i}(\Delta_i)]}{q^\alpha} -\sum_{i=1}^n \E\left[g_i(\Theta_i,\,\Delta_i)\,\overline{\mu}\left( \dfrac {p}{\Theta_i},\;\dfrac {q}{\Delta_i}\right) \right] \Bigg)\,.
\eeam
\end{enumerate}
\ethe
	
\pr~
\begin{enumerate}
\item[(i)]
Relation~\eqref{eq.AKP.6.14} implies that
	\beao
\PP[S_n^{\Theta}> p\,x\,,\; T_n^{\Delta}>q\,x] \leq \psi_{and}(x,\,p,\,q;\,n) \leq \PP[S_n^{\Theta+}> p\,x\,,\; T_n^{\Delta+}>q\,x]\,,
	\eeao
	which together with Lemma \ref{lem.AKP.5.1} gives  \eqref{eq.AKP.6.15}, for $\psi_{and}(x,\,p,\,q;\,n) $.

Further we deal with $\psi_{sim}(x,\,p,\,q;\,n)$. The lower bound result is derived as follows
\beao
\psi_{sim}(x,\,p,\,q;\,n) &=& \PP\left[\max_{1\leq k \leq n} \left( \sum_{i=1}^k \Theta_{i}\,X_i >p\,x\,,\; \sum_{i=1}^k \Delta_{i}\,Y_i > q\,x \right) \right]\\[2mm] \notag
&\geq& \PP\left[\sum_{i=1}^n \Theta_{i}\,X_i > p\,x\,,\;  \sum_{i=1}^n \Delta_{i}\,Y_i > q\,x\right] \\[2mm] \notag
&\sim& \bV(x)\,\sum_{i=1}^n \E\left[g_i(\Theta_i,\,\Delta_i)\,\overline{\mu}\left( \dfrac {p}{\Theta_i},\;\dfrac {q}{\Delta_i}\right) \right]\,, 
\eeao
where at the last step we use Lemma~\ref{lem.AKP.5.1}. 

The upper bound result is true, due to \eqref{eq.AKP.6.15} for $\psi_{and}(x,\,p,\,q;\,n)$ and the fact that $\psi_{sim}(x,\,p,\,q;\,n) \leq \psi_{and}(x,\,p,\,q;\,n)$ for any triple $(x,\,p,\,q)$.
	
\item[(ii)]
The upper asymptotic bound for $\psi_{or}(x,\,p,\,q;\,n)$ is true given that
\beam \label{eq.AKP.6.19} \notag
&&\psi_{or}(x,\,p,\,q;\,n)= \PP\Bigg[\left\{\max_{1\leq k \leq n} \sum_{i=1}^k \Theta_{i}\,X_i>  p\,x \right\} \bigcup \left\{\max_{1\leq l \leq n} \sum_{i=1}^l \Delta_i\,Y_i> q\,x\right\}\Bigg]\\[2mm] 
&&\leq  \PP\left[\left\{ \sum_{i=1}^n \Theta_{i}\,X_i^+>  p\,x \right\} \bigcup \left\{ \sum_{i=1}^n \Delta_i\,Y_i^+> q\,x\right\}\right]\\[2mm] \notag
&&= \PP\left[ \sum_{i=1}^n \Theta_{i} X_i^+> p x \right] + \PP\left[ \sum_{i=1}^n \Delta_{i} Y_i^+> q x \right]-\!\PP\left[\sum_{i=1}^n \Theta_{i} X_i^+> p x ,\,\sum_{i=1}^n \Delta_{i} Y_i^+> q x\right]  \\[2mm] \notag
&&\sim \bV(x)\,\Bigg( \sum_{i=1}^n\dfrac{\E[\Theta_i^{\alpha}\,h_{1,i}(\Theta_i)]}{p^\alpha}  + \sum_{i=1}^n\dfrac{\E[\Delta_i^{\alpha}\,h_{2,i}(\Delta_i)]}{q^\alpha}-\sum_{i=1}^n \E\left[g_i(\Theta_i,\,\Delta_i)\,\overline{\mu}\left( \dfrac {p}{\Theta_i},\;\dfrac {q}{\Delta_i}\right) \right] \Bigg)\,
\eeam
where we use \eqref{eq.AKP.4.42} and Lemma \ref{lem.AKP.5.1} at the last step, since we consider here the standard $BRV$. 

For the lower asymptotic bound of $\psi_{or}(x,\,p,\,q;\,n)$ by  \eqref{eq.AKP.4.42} and Lemma \ref{lem.AKP.5.1} we find that 
\beam \label{eq.AKP.6.20}
&&\psi_{or}(x,\,p,\,q;\,n) \geq \PP\Bigg[\left\{\sum_{i=1}^n \Theta_{i}\,X_i>  p\,x \right\} \bigcup \left\{ \sum_{i=1}^n \Delta_i\,Y_i> q\,x\right\}\Bigg]\\[2mm] \notag
&&=\PP\left[\sum_{i=1}^n \Theta_{i}\,X_i> p\,x \right]+ \PP\left[ \sum_{i=1}^n \Delta_{i}\,Y_i> q\,x \right]-\PP\left[\sum_{i=1}^n \Theta_{i}\,X_i>  p\,x \,,\; \sum_{i=1}^n \Delta_{i}\,Y_i>   q\,x\right]  \\[2mm] \notag
&&\sim \bV(x)\,\Bigg( \sum_{i=1}^n\dfrac{\E[\Theta_i^{\alpha}\,h_{1,i}(\Theta_i)]}{p^\alpha}  + \sum_{i=1}^n\dfrac{\E[\Delta_i^{\alpha}\,h_{2,i}(\Delta_i)]}{q^\alpha} -\sum_{i=1}^n \E\left[g_i(\Theta_i,\,\Delta_i)\,\overline{\mu}\left( \dfrac {p}{\Theta_i},\;\dfrac {q}{\Delta_i}\right) \right] \Bigg)\,.
\eeam
From relations \eqref{eq.AKP.6.19} and \eqref{eq.AKP.6.20} we get relation \eqref{eq.AKP.6.17}.~\halmos
\end{enumerate}

By Theorem \ref{th.AKP.5.1} we can see that the ruin probabilities $\psi_{and}(x,\,p,\,q;\,n)$ and $\psi_{sim}(x,\,p,\,q;\,n)$ are asymptotically equivalent, that indicates that their difference is negligible, that means the simultaneous and non-simultaneous ruin  in both portfolios under the presence of heavy tails on the net losses, have almost the same probability. This fact can be observed also in several papers on continuous time risk models on more dimensions, under standard $MRV$ condition, see for example \cite{cheng:konstantinides:wang:2022}, \cite{cheng:konstantinides:wang:2024}, \cite{yang:su:2023}. Further we recall that relations \eqref{eq.AKP.6.12} and \eqref{eq.AKP.6.13} concept immediately the impact of the dependence among the $X_i,\,Y_i,\,\Theta_i,\,\Delta_i$ via the functions $h_{1,i},\,h_{2,i},\,g_{i}$.

\subsection{Joint Expected Shortfall} \label{subsec.KP.6.2}

In this section we study an application of our main results for to quantify the marginal effect of concomitant extreme events. This tail measure is inspired by the so-called \emph{Joint Expected Shortfall (JES)} risk measure, introduced in \cite{ji:tan:yang:2021} and given as follows:
\begin{eqnarray}\label{JES_def}
JES(x):=\E\left[\Theta X\;\big|\;\Theta X>U_F(x),\,\Delta Y> U_G(x)\right]
\end{eqnarray}
and we evaluate $JES(x)$ as $\xto$.  The main result of this section is given in Proposition~\ref{JES_prop}. Let notice that the limit measure $\widehat{\mu}$ in the following proposition was given in \eqref{eq.AKP.2.19}. We assume that the $(\Theta\,X,\,\Delta\,Y)$ are asymptotic dependent through the condition \eqref{eq.AKP.6.0}. We note that in some cases (via the last term of \eqref{eq.AKP.2.19}), such an assumption on $\widehat{\mu}$, is implied by the condition $\overline{\mu}(1,\,1)> 0$.

\begin{proposition} \label{JES_prop}
Let $(X,\,Y)$ and $(\Theta,\,\Delta)$ be nonnegative random vectors such that $(X,\,Y) \in BRV_{-\alpha,\,-\beta}(\mu)$ holds for some finite $\alpha, \,\beta >0$, and assume that all the conditions in Proposition~\ref{th.AKP.3.1} hold, with $\inf_{\theta \in S_{\Theta},\,\delta \in S_{\Delta}}g(\theta,\,\delta) >0 $.  Moreover, assume that there exists $x_0>0$ such that 
\beao
JES(x)<\infty\,,
\eeao 
for any $x>x_0$.  Then, for any $\alpha>1$ we have that 
\beao
JES(x) \sim U_F(x)\left[
1+\int_1^\infty \dfrac{ \E\left[g(\Theta,\,\Delta)\,\overline{\mu}\Big(\dfrac{\xi}{\Theta},\,\dfrac{1}{\Delta}\Big)\right]}{ \E\left[g(\Theta,\,\Delta)\,\overline{\mu}\Big(\dfrac{1}{\Theta},\,\dfrac{1}{\Delta}\Big)\right]}\,d\xi \right],
\eeao
provided that $\widehat{\mu}\big((1,\infty]\times(1,\infty]\big)>0$.
\end{proposition}

\pr~ 
Since $JES(x)<\infty$ for any $x>x_0$, one may find that
\beao
\dfrac{JES(x)}{U_F(x)}&=&\frac{1}{U_F(x)}\int_0^\infty \Pr\big(\Theta X>t|\Theta X>U_F(x),\,\Delta Y> U_G(x)\big)\;dt \nonumber\\ \notag
&=& 1+\frac{1}{U_F(x)}\int_{U_F(x)}^\infty \frac{\Pr\big(\Theta X>t,\;\Delta Y> U_G(x)\big)}{\Pr\big(\Theta X>U_F(x),\,\Delta Y> U_G(x)\big)}\;dt\\ \notag
&=& 1 +\int_1^\infty \frac{\Pr\big(\Theta X>U_F(x)\xi,\;\Delta Y> U_G(x)\big)}{\Pr\big(\Theta X>U_F(x),\,\Delta Y> U_G(x)\big)}\;d{\xi} \nonumber\\
&\sim& 1+\int_1^\infty \dfrac{ \E\left[g(\Theta,\,\Delta)\,\overline{\mu}\Big(\dfrac{\xi}{\Theta},\,\dfrac{1}{\Delta}\Big)\right]}{ \E\left[g(\Theta,\,\Delta)\,\overline{\mu}\Big(\dfrac{1}{\Theta},\,\dfrac{1}{\Delta}\Big)\right]}\;d\xi\,,
\eeao
where a change of variables, $t=\xi U_F(x)$, is made at the third step, while the Dominated Convergence Theorem and the relation~\eqref{eq.AKP.2.25}  are applied at the last step.  The latter follows from the relation
\beao \label{JES_def_eq2}
&&\int_1^\infty \dfrac{\Pr\big(\Theta X>U_F(x)\xi,\;\Delta Y> U_G(x)\big)}{\Pr\big(\Theta X>U_F(x),\,\Delta Y> U_G(x)\big)}\;d\xi \\[2mm] \notag
&&\leq \dfrac{\Pr\big(\Theta X>U_F(x)\big)}{\Pr\big(\Theta X>U_F(x),\,\Delta Y> U_G(x)\big)}\int_1^\infty \frac{\Pr\big(\Theta X>U_F(x)\xi\big)}{\Pr\big(\Theta X>U_F(x)\big)}\;d\xi \\[2mm] \notag
&&\leq K\,\frac{\widehat{\mu}\big((1,\infty]\times [0,\infty]\big)}{\widehat{\mu}\big((1,\infty]\times(1,\infty]\big)}\,\int_1^\infty \xi^{\epsilon_2-\alpha}\;d\xi,\nonumber
\eeao
where the last step follows from Theorem~\ref{th.AKP.3.1} and the \cite[Lem. 1 (b)]{li:2018} so that there exists $x_0>0$ such that last inequality holds for some $K>0$ and some $0<\epsilon_2<\alpha - 1$. These enable to apply the Dominated Convergence Theorem, which completes this proof.
~\halmos

\noindent \textbf{Acknowledgments.} 
We would like to thank Vali Asimit, for his helpful comments, that improved substantially the text.

\noindent \textbf{Funding Declaration.} 
No competing or funding interests, that influence the results of this paper.

\end{document}